\documentclass{amsart}
\usepackage{amsfonts,amssymb,amsmath,amsthm}
\usepackage{mathrsfs}
\usepackage[all]{xy}
\usepackage{url}
\usepackage{enumerate}

\newtheorem{thrm}{Theorem}[section]

\newtheorem{prop}[thrm]{Proposition}
\newtheorem{cor}[thrm]{Corollary}
\theoremstyle{definition}
\newtheorem{definition}[thrm]{Definition}
\newtheorem{remark}[thrm]{Remark}
\newtheorem{example}[thrm]{Example}
\numberwithin{equation}{section}
\author{L. Vitagliano}
\address{
DipMat, University of Salerno\\
and Istituto Nazionale di Fisica Nucleare, GC Salerno\\
via Ponte don Melillo\\
84084 Fisciano (SA) Italy}
\email{lvitagliano@unisa.it}

\begin{document}

\title[Partial Differential Hamiltonian Systems]{Partial Differential Hamiltonian Systems}

\begin{abstract}
We define partial differential (PD in the following), i.e., field theoretic analogues of Hamiltonian systems on abstract symplectic manifolds and study their main properties, namely, PD Hamilton equations, PD Noether theorem, PD Poisson bracket, etc.. Unlike in standard multisymplectic approach to Hamiltonian field theory, in our formalism, the geometric structure (kinematics) and the dynamical information on the \textquotedblleft phase space\textquotedblright\ appear as just different components of one single geometric object.
\end{abstract}
\maketitle

\section{Introduction}

First order Lagrangian mechanics can be naturally generalized to higher
order Lagrangian field theory. Moreover, the latter can be presented in a very
elegant and precise algebro-geometric fashion \cite{v84}.
In particular, it is clear what all the involved geometric structures
(higher order jets, Cartan distribution, $\mathscr{C}$-spectral sequence, etc.,
\cite{v84,b...99}) are. On the other hand it seems to be quite hard to understand
what the most \textquotedblleft reasonable, unambiguous, higher order, field
theoretic generalization\textquotedblright\ of Hamiltonian mechanics on
abstract symplectic manifolds is. Actually, there exists a universally accepted
generalization of the standard mechanical picture
\[
\text{Lagrangian mechanics on }TQ\Longrightarrow\text{ Hamiltonian
mechanics on }T^{\ast}Q\text{,}%
\]
$Q$ being a smooth manifold, to the picture
\begin{equation}
\text{Lagrangian field theory on }J^{1}\pi\Longrightarrow\text{ Hamiltonian
field theory on }\mathscr{M}\pi\text{,} \label{Picture}%
\end{equation}
$\pi$ being a fiber bundle, $J^{1}\pi$ its first jet space and $\mathscr{M}\pi
$ its multimomentum space \cite{gs73} (see also \cite{hk04}, \cite{r05} for a recent
review, and \cite{g12} for an approach \textit{\`{a} l\`{a}} Tulczyjew). Picture (\ref{Picture}) includes, in particular, a generalization of
the Legendre transform. Along this path an analogous
structure of the symplectic structure on $T^{\ast}Q$, which has been named the
\emph{multisymplectic structure of} $\mathscr{M}\pi$ (see, for instance,
\cite{gim98}), has been discovered. A whole literature exists about properties of such
structure, which is generically referred to as \emph{multisymplectic
geometry} of $\mathscr{M}\pi$ (see references in \cite{r05}). In particular,
efforts were made to find multisymplectic analogues of all properties of
$T^{\ast}Q$ (including, for instance, the Poisson bracket
\cite{ks75,k97,fr01,fpr03,fpr03b}). Now, it is natural to wonder if it is
possible to reasonably further generalize in two different directions. The
first one is towards a picture
\begin{equation}
\text{Lagrangian field theory on }J^{\infty}\pi\Longrightarrow\text{ higher
order Hamiltonian field theory,} \label{Picture2}%
\end{equation}
$J^{\infty}\pi$ being the $\infty$th jet space of $\pi$, including a higher
order generalization of the Legendre transform. There is no universally
accepted answer about picture (\ref{Picture2})
(see, for instance, \cite{d77,aa82,s82,k84b,sc90,s91,k02,av04} and references therein). Most
often they involve the choice of some extra structure other than the natural
ones on $J^{\infty}\pi$. Recently, we proposed in \cite{v09b} an answer that is free from such ambiguities. 

The second direction in which to generalize picture (\ref{Picture}) can be
illustrated as follows. $T^{\ast}Q$ is just a very special example of
(pre)symplectic manifold. Actually Hamiltonian mechanics can (and should,
in some cases \cite{gnh78}) be formulated on abstract (pre)symplectic manifolds.
Similarly, it is natural to wonder if there exists the concept of abstract
multi(pre)symplectic manifolds in such a way that Hamiltonian field theory
could be reasonably formulated on them. In the literature there can be found some
proposals of should be abstract multi(pre)symplectic manifolds (see, for
instance, \cite{a92b,cid99}). In particular, definitions have been given in
such a way to be able to prove multisymplectic analogues of the celebrated
Darboux lemma \cite{m88,fg08}. The recent definitions by Forger and Gomes
\cite{fg08} appear to be the most satisfactory in that they are
\textquotedblleft minimal\textquotedblright\ on one side and duly model in an
abstract fashion the relevant geometric properties of $\mathscr{M}\pi$ on the other side. In
their work Forger and Gomes illustrate, in particular, the role played by
fiber bundles in the should be definition of multi(pre)symplectic structure.
The next step forward should be to formulate Hamiltonian field theory on
multisymplectic bundles.

In this paper we present our own proposal about what should be an abstract,
first order, Hamiltonian field theory. We call such proposal the \emph{theory
of partial differential} (\emph{PD} in the following)
\emph{Hamiltonian systems} so to 1) stress that it is a natural generalization
of the theory of Hamiltonian systems on abstract symplectic manifolds, 2)
distinguish it from the special case of Hamiltonian field theory on
$\mathscr{M}\pi$. A PD Hamiltonian system encompasses both the kinematics
(encoded, in picture (\ref{Picture}), by the multisymplectic structure in
$\mathscr{M}\pi$) and the dynamics (encoded, in picture (\ref{Picture}), by the
so called \emph{Hamiltonian section }\cite{r05}) which appear as just
different components of one single geometric object. Namely, the main
difference between a PD Hamiltonian system and a multi(pre)symplectic
structure (whatever the reader understand for this) is the dynamical content
of the former (as opposed to the just kinematical one of the latter). Notice that this idea is already present in literature \cite{k73}. However, our formalism differs from the one in \cite{k73} in that it is adapted to the fibered structure of the manifold of \textquotedblleft field variables\textquotedblright .

As already mentioned, standard examples of PD Hamiltonian systems come from
Lagrangian field theory. Consider a field theory on a \textquotedblleft
space-time\textquotedblright\ $M$ with coordinates $x:=(\ldots,x^{i},\ldots)$,
defined by a Lagrangian density
\[
\mathscr{L}=L(x,u,u^{\prime})d^{n}x
\]
depending on some field variables $u:=(\ldots,u^{\alpha},\ldots)$ and their
partial derivatives $u^{\prime}:=(\ldots,u_{i}^{\alpha},\ldots)$.
The Lagrangian density $\mathscr{L}$ determines a Legendre transform $F\mathscr{L}$ (see, for
instance, \cite{r05} for a geometric formulation of the constructions in this
paragraph) from the space $J^{1}$ of the $(x,u,u^{\prime})$'s to the so called
\emph{multimomentum space} $J^{\dag}$ with coordinates $(x,u,p)$,
$p:=(\ldots,p_{\alpha}^{i},\ldots)$ being the \emph{multimomenta}.
The Legendre transform $F\mathscr{L}:J^{1}\longrightarrow J^{\dag}$ is defined as
\[
F\mathscr{L}(x,u,u^{\prime}):=(x,u,\partial L/\/\partial u^{\prime}).
\]
Moreover, $\mathscr{L}$ determines in a canonical way the following
$(n+1)$-form on $J^{1}$:
\[
\omega_{\mathscr{L}}:=d\tfrac{\partial L}{\partial u_{i}^{\alpha}}\wedge
du^{\alpha}\wedge d^{n-1}x_{i}-dE\wedge d^{n}x,\quad E:=u_{i}^{\alpha}%
\tfrac{\partial L}{\partial u_{i}^{\alpha}}-L.
\]
The form $\omega_{\mathscr{L}}$ is \textquotedblleft constant\textquotedblright\ along the fibers of $F\mathscr{L}$ and,
therefore, determines an $(n+1)$-form $\omega$ on the submanifold
$C_{0}:=\operatorname{im}F\mathscr{L}\subset J^{\dag}$. The forms $\omega_{\mathscr{L}}$
and $\omega$ are standard examples of PD Hamiltonian systems. In a way that
will be clear later, they determine PDEs for sections of the bundles
$J^{1}\longrightarrow M$ and $C_{0}\longrightarrow M$. When $F\mathscr{L}$ is
invertible,
\[
\omega = dp_{\alpha}^{i}\wedge du^{\alpha}\wedge d^{n-1}x_{i}-dH\wedge
d^{n}x,\quad H:=E\circ FL^{-1}%
\]
and the above mentioned PDEs are the Euler-Lagrange and the de Donder-Weyl
equations of the theory, respectively. Using the theory of PD Hamiltonian systems, one can generalize these considerations to more general field theories, in particular those depending on higher derivatives of the fields \cite{v09b}.

The paper is divided into 8 sections. In Section \ref{SecNotations} we
collect our notations and conventions and recall basic differential geometric
facts that will be used in the main part of the paper.

In Sections \ref{SecAffForms} we define what we call \emph{affine forms} on
fiber bundles. The introduction of affine forms can be motivated as follows.
Trajectories in Hamiltonian mechanics are curves, whose 1st derivative at a
point is naturally understood as a tangent vector. In their turn, tangent
vectors can be inserted into differential forms and, in particular, a
symplectic one, and Hamilton equations are written in terms of such an insertion. Trajectories in field theory are sections of a fiber bundle
$\alpha:P\longrightarrow M$, whose 1st derivative at a point is naturally
understood as a point in $J^{1}\alpha$. In their turn, points of $J^{1}\alpha$
can be inserted into affine forms and, in particular, a PD Hamiltonian system
(see Section \ref{SecPDHamSys}), and PD Hamilton equations are written in terms of such an insertion. Recall now that the natural projection
$J^{1}\alpha\longrightarrow P$ is an affine bundle whose sections are
naturally interpreted as (Ehresmann) connections in $\alpha$. Thus, connections
and affine geometry play a prominent role in the theory of PD Hamiltonian
systems. The affine geometry is hidden in standard Hamiltonian mechanics by an
\textit{a priori} choice of the parametrization of the time axis (see
\cite{ggu04,ggu06}, and references therein, for the role of affine geometry in theoretical
mechanics). Similarly, even if the role of connections in field theory has
been often recognized (see, for instance, \cite{s87,emr00}), their affine
geometry is sometimes hidden in Hamiltonian field theory on $\mathscr{M}\pi$
by the use of multivectors, or even decomposable ones \cite{pr02,pr02b} (which
is just a multidimensional analogue of choosing a parameterization of time).
Actually, we show in Subsection \ref{SecAff} that affine forms can be
understood as standard differential forms of a special kind. Nevertheless, we prefer to keep the distinction for
foundational reasons.

In Section \ref{SecCalculus} we discuss standard operations with affine forms.
Essentially because of the interpretation of affine forms as standard
differential forms we mentioned above, some of this operations (for instance,
the insertion of a connection into an affine form \cite{emr96}, or the
differential of an affine form) were actually already defined in the
literature, or can be understood as standard operations with forms. We stress
again that we will keep the distinction. Finally, we discuss relevant affine
form cohomologies proving an affine form version of the Poincar\'{e} lemma.

In Section \ref{SecPDHamSys} we introduce PD (pre)Hamiltonian systems, discuss
their geometry and the geometry of the associated PD Hamilton equations, with
some references to the singular, constrained case (see
\cite{dmm96,dmm96b,d...02,d...05} for an account of the constraint algorithm
in first order field theory). For completeness, we also relate PD Hamiltonian
systems to multi(pre)symplectic structures \textit{\`{a} l\`{a} }Forger
\cite{fg08} and the calculus of variations. Notice that the theory of PD Hamiltonian systems
is somehow in between (abstract) multisymplectic field theory and polysymplectic field theory.

In Section \ref{SecNoether} we introduce PD Noether symmetries and currents of
a PD Hamiltonian system. In view of the dynamical content of the latter we are
able to prove a Noether theorem (see also \cite{dms04}). Moreover, there is a natural Lie bracket
(named PD Poisson bracket) among PD Noether currents. As already mentioned,
multisymplectic analogues of the Poisson bracket have been already discussed in
the literature \cite{ks75,k97,fr01,fpr03,fpr03b}. However, we emphasize here the 
dynamical nature of PD Poisson bracket (see also \cite{fr05,v09}). Namely, such bracket is just part
of the Peierls bracket \cite{v09} among conservation laws of the underlying
Lagrangian theory and we don't try to extend it to non-conserved currents.
Indeed, our opinion is that the existence of a Poisson bracket among
non-conserved functions in Hamiltonian mechanics is essentially due to the
existence of a preferred Hamiltonian system on any symplectic manifold $N$,
i.e., the one with $0$ Hamiltonian, for which every function on $N$ is a
conservation law. Finally we discuss the (gauge) reduction of a degenerate
(but unconstrained) PD Hamiltonian system.

In Section \ref{SecComp} we propose few examples of PD Hamiltonian systems,
including the computation of their PD Noether symmetries and currents or,
in one case, their reduction.

We conclude with Section \ref{SecMPG} where we briefly discuss the emergence of PD Hamiltonian systems in Mathematical Physics and Geometry, providing further motivations for their introduction.

\section{Notations and Conventions\label{SecNotations}}

In this section we collect notations and conventions about some general
constructions that will be used in the following.

Let $N$ be a manifold. We denote by $C^{\infty}(N)$ the $\mathbb{R}$-algebra
of smooth, $\mathbb{R}$-valued functions on $N$. A vector field $X$ over $N$
will be always understood as a derivation $X:C^{\infty}(N)\longrightarrow
C^{\infty}(N)$. We denote by $\mathrm{D}(N)$ the $C^{\infty}(N)$-module of
vector fields over $N$, by $\Lambda(M)=\bigoplus_{k}\Lambda^{k}(N)$ the graded
$\mathbb{R}$-algebra of differential forms over $N$, by $d:\Lambda
(N)\longrightarrow\Lambda(N)$ the de Rham differential, and
by $H(N)=\bigoplus_{k}H^{k}(N)$ the de Rham cohomology. If $F:N_{1}%
\longrightarrow N$ is a smooth map of manifolds, we denote by $F^{\ast
}:\Lambda(N)\longrightarrow\Lambda(N_{1})$ its pull-back. We will understand
everywhere the wedge product $\wedge$ of differential forms, i.e., for
$\omega,\omega_{1}\in\Lambda(N)$, instead of writing $\omega\wedge\omega_{1}$,
we will simply write $\omega\omega_{1}$. We assume the reader to be familiar
with Fr\"{o}licher-Nijenhuis calculus on form valued vector fields (insertion
$i_{Z}\omega$ of a form valued vector field $Z$ into a differential form
$\omega$, Lie derivative $L_{Z}\omega$ of a differential form $\omega$ along a
form valued vector fields $Z$, Fr\"{o}licher-Nijenhuis bracket, etc., see, for instance, \cite{m08}).

Let $\varpi:W\longrightarrow N$ be an affine bundle (or, possibly, a vector
bundle) and let $F:N_{1}\longrightarrow N$ be a smooth map of manifolds. The affine
space of smooth sections of $\varpi$ will be denoted by $\Gamma(\varpi)$. For
$x\in N$, we put, sometimes, $\Gamma(\varpi)|_{x}:=\varpi^{-1}(x)$ and, for
$\chi\in\Gamma(\varpi)$, we also put $\chi_{x}:=\chi(x)$. The affine bundle on
$N_{1}$ induced by $\varpi$ via $F$ will be denote by $\varpi|_{F}%
:W|_{F}\longrightarrow N$:
\[%
\begin{array}
[c]{c}%
\xymatrix{W|_F \ar[r] \ar[d]_-{\varpi|_F} & W \ar[d]^-{\varpi} \\ N_1 \ar[r]^-F & N }\end{array}
.
\]
We also denote $\Gamma(\varpi)|_F := \Gamma(\varpi|_F)$. For any section $s\in\Gamma(\varpi)$, there exists a unique section, which,
abusing the notation, we denote by $s|_{F}\in\Gamma(\varpi|_{F})$, such that
the diagram
\[
\xymatrix{W|_F \ar[r]  & W  \\
N_1 \ar[r]^-F    \ar[u]^-{s|_F}                    &  N \ar[u]_-{s}}
\]
commutes. Elements in $\Gamma(\varpi)|_{F}$ are called \emph{sections of
}$\varpi$\emph{ along }$F$. If $F$ is an embedding $\varpi|_{F}$,
$\Gamma(\varpi)|_{F}$ and $s|_{F}$ will be referred to as \emph{the
restriction to }$N_{1}$ of $\varpi,$ $\Gamma(\varpi)$ and $s$, respectively.
If $\varpi_{1}:W_{1}\longrightarrow N$ is an other affine bundle and
$A:\Gamma(\varpi)\longrightarrow\Gamma(\varpi_{1})$ is an affine map then
there exists a unique affine map $A|_{F}:\Gamma(\varpi)|_{F}\longrightarrow
\Gamma(\varpi_{1})|_{F}$ such that $A|_{F}(s|_{F})=A(s)|_{F}$ for all
$s\in\Gamma(\varpi)$.

Let $\alpha:P\longrightarrow M$ be a fiber bundle. A vector field
$X\in\mathrm{D}(P)$ is called $\alpha$\emph{-projectable} iff there exists
$\check{X}\in\mathrm{D}(P)$ such that $X\circ\alpha^{\ast}=\alpha^{\ast}%
\circ\check{X}$. The vector field $\check{X}$ is called the $\alpha$\emph{-projection of
}$X$. Vector fields that are $\alpha$-projectable form a Lie subalgebra in
$\mathrm{D}(P)$ denoted by $\mathrm{D}_{V}(P,\alpha)$ (or simply
$\mathrm{D}_{V}$ if this does not lead to confusion). An $\alpha$-projectable
vector field projecting onto the $0$ vector field is an $\alpha$\emph{-vertical
vector field}. Vector fields that are $\alpha$-vertical form an ideal in $\mathrm{D}%
_{V}$ denoted by $V\mathrm{D}(P,\alpha)$ (or simply $V\mathrm{D}$). Notice
that, if $\alpha$ has connected fiber, then $\mathrm{D}_{V}$ is the stabilizer
of $V\mathrm{D}$ in $\mathrm{D}(P)$, i.e., $\mathrm{D}_{V}=\{X\in
\mathrm{D}(P)\;|\;[X,V\mathrm{D}]\subset V\mathrm{D}\}$.

Let $\alpha:P\longrightarrow M$ be as above, $\dim M=n$, $\dim P=m+n$. Denote
by $\alpha_{1}:J^{1}\alpha\longrightarrow M$ the bundle of $1$-jets of local
sections of $\alpha$ \cite{s89,b...99}, and by $\alpha_{1,0}:J^{1}\alpha
\longrightarrow P$ the canonical projection. For any local section
$\sigma:U\longrightarrow P$ of $\alpha$, $U\subset M$ being an open subset, we
denote by $\dot{\sigma}:U\longrightarrow J^{1}\alpha$ its $1$st jet
prolongation. Any system of adapted to $\alpha$ coordinates $(\ldots
,x^{i},\ldots,y^{a},\ldots)$ on $P$, $\ldots,x^i,\ldots$ being coordinates on $M$
and $\ldots,y^a,\ldots$ fiber coordinates on $P$, gives rise to the system 
of jet coordinates $(\ldots,x^{i},\ldots,u^{a},\ldots,y_{i}^{a},\ldots)$ on
$J^{1}\alpha$, $i=1,\ldots,n$, $a=1,\ldots,m$. Recall that $\alpha_{1,0}$ is
an affine bundle and a section $\nabla:P\longrightarrow J^{1}\alpha$ of it is
naturally interpreted as a (Ehresmann) connection in $\alpha$. We assume the
reader to be familiar with the geometry of connections (see, for instance,
\cite{m08}). A connection $\nabla$ is locally represented as
\[
\nabla:y_{i}^{a}=\nabla_{i}^{a},
\]
$\ldots,\nabla_{i}^{a},\ldots$ being local functions on $P$. The space $\Gamma(\alpha_{1,0})$ of all such
sections will be also denoted by $C(P,\alpha)$ (or simply $C$).

Let $\alpha:P\longrightarrow M$ be as above, $\alpha^{\prime}:P^{\prime
}\longrightarrow M$ be another fiber bundle and let $G:P\longrightarrow P^{\prime}$ be a
bundle morphism (over the identity $\operatorname*{id}_{M}:M\longrightarrow
M$), i.e., a smooth map such that $\alpha^{\prime}\circ G=\alpha$. First of
all, recall that there exists a unique bundle morphism $j_{1}G:J^{1}%
\alpha\longrightarrow J^{1}\alpha^{\prime}$ such that $j_{1}G\circ\dot{\sigma
}=(G\circ\sigma)\dot{\,}$ for all local sections $\sigma$ of $\alpha$.
The bundle morphism $j_{1}G$ is the \emph{first jet prolongation of }$G$ and diagram
\[
\xymatrix{ J^1 \alpha \ar[rr]^-{j_1 G} \ar[d]_-{\alpha_{1,0}} &   & J^1 \alpha^\prime  \ar[d]^-{\alpha^\prime_{1,0}} \\
              P       \ar[rr]^-{G}     \ar[rd]_-{\alpha}      &   &      P^\prime      \ar[dl]^-{\alpha^\prime} \\
                                                              & M & }
\]
commutes. Now, a connection $\nabla\in C(P,\alpha)$ and a connection
$\nabla^{\prime}\in C(P^{\prime},\alpha^{\prime})$ are said $G$%
\emph{-compatible} iff $\nabla^{\prime}\circ G=j_{1}G\circ\nabla$.

Let
\[
\xymatrix{\cdots \ar[r] & K_{l-1} \ar[r]^-{\delta_{l-1}} & K_{l} \ar[r]^-{\delta_{l}} & K_{l+1} \ar[r]^-{\delta_{l+1}} & \cdots}
\]
be a complex. Put $K:=\bigoplus_{l}K_{l}$ and $\delta:=\bigoplus_{l}\delta
_{l}$. We denote by $H(K,\delta):=\bigoplus_{l}H^{l}(K,\delta)$, the
cohomology space of $(K,\delta)$, $H^{l}(K,\delta):=\ker\delta_{l}%
/\operatorname{im}\delta_{l-1}$.

Let $A$ be a commutative $\mathbb{R}$-algebra, $\boldsymbol{M},\boldsymbol{M}%
_{1}$ be $A$-modules and let $\boldsymbol{A}$ be an affine space modeled over
$\boldsymbol{M}$. We denote by $\mathrm{Aff}_{A}(\boldsymbol{A},\boldsymbol{M}%
_{1})$ (resp.{} $\mathrm{Hom}_{A}(\boldsymbol{M},\boldsymbol{M}_{1})$) the
$A$-module of affine (resp.{} $A$-linear) maps $\boldsymbol{A}\longrightarrow
\boldsymbol{M}_{1}$ (resp.{} $\boldsymbol{M}\longrightarrow\boldsymbol{M}_{1}$).
If $\phi\in\mathrm{Aff}_{A}(\boldsymbol{A},\boldsymbol{M}_{1})$, its
\emph{linear part} $\underline{\phi}$ is an element in $\mathrm{Hom}%
_{A}(\boldsymbol{M},\boldsymbol{M}_{1})$.

Let $m,r$ be positive integers and let $\ldots,A_{a_{1}\cdots a_{r}},\ldots$ be
elements in a real vector space, $a_{1},\ldots,a_{r}=1,\ldots,m$. We denote by
$\ldots,A_{[a_{1}\cdots a_{r}]},\ldots$ their skew-symmetrization, i.e.,
\[
A_{[a_{1}\cdots a_{r}]}:=\tfrac{1}{s!}\sum_{\sigma\in S_{r}}\varepsilon
(\sigma)A_{a_{\sigma(1)}\cdots a_{\sigma(r)}},
\]
$S_{r}$ being the group of permutations of $\{1,\ldots,r\}$ and $\varepsilon
(\sigma)$ the sign of $\sigma\in S_{r}$.

We denote by $\simeq$ (resp.{} $\approx$) a canonical (resp.{} non-canonical)
isomorphism between algebraic structures and by $\equiv$ an equivalence of
notations. For instance, for $\alpha:P\longrightarrow M$ as above, $V\mathrm{D} \equiv V\mathrm{D}(P,M)$. Finally, we understand the sum over upper-lower pairs of repeated indexes.

\section{Affine Forms on Fiber Bundles\label{SecAffForms}}

\subsection{Special Forms on Fiber Bundles}

Let $\alpha:P\longrightarrow M$ be a fiber bundle, $A:=C^{\infty}(P)$,
$A_{0}:=C^{\infty}(M)$, $x^{1},\ldots,x^{n}$ be coordinates on $M$, $\dim M=n$,
and let $y^{1},\ldots,y^{m}$ be fiber coordinates on $P$, $\dim P=n+m$. In the following
we will often understand the monomorphism of algebras $\alpha^\ast:A_0 \longrightarrow A$, whose image is made of functions on $P$ which are constant along the fibers of $\alpha$.
The space $\mathrm{D}_{V}$ (resp.{} $V\mathrm{D}$) is made of vector
fields $X$ locally of the form $X=X^{i}\partial_{i}+Y^{a}\partial_{a}$ (resp.
$X=Y^{a}\partial_{a}$) where $X^{i}=X^{i}(x^{1},\ldots,x^{n})$, $\partial
_{i}:=\partial/\partial x^{i}$, $i=1,\ldots,n$, $\partial_{a}=\partial
/\partial y^{a}$, $a=1,\ldots,m$.

Denote by $\Lambda_{1}(P,\alpha)=\bigoplus_{k}\Lambda_{1}^{k}(P,\alpha)$ (or
simply $\Lambda_{1}=\bigoplus_{k}\Lambda_{1}^{k}$) the differential (graded)
ideal of differential forms on $P$ vanishing when pulled-back to fibers of
$\alpha$, i.e., $\omega \in \Lambda^k_1$, $k \geq 0$ iff $\omega \in \Lambda^k(P)$ and 
$i^\ast_{\alpha^{-1}(x)}(\omega)=0$ for all $x\in M$, $i_{\alpha^{-1}(x)}: 
\alpha^{-1}(x) \longrightarrow P$ being the embedding of the fiber $\alpha^{-1}(x)$
of $\alpha$ through $x \in M$. Moreover, denote by $\Lambda_{p}(P,a)=\bigoplus_{k}\Lambda_{p}^{k}(P,\alpha)$
(or simply $\Lambda_{p}=\bigoplus_{k}\Lambda_{p}^{k}$) the $p$-th exterior
power of $\Lambda_1$. For all $k$ and $p$, $\Lambda_{p}^{k}$ is made of differential $k$-forms $\omega$ such that
$(i_{Y_{1}}\circ\cdots\circ i_{Y_{k-p+1}})\omega=0$ for every $Y_{1}%
,\ldots,Y_{k-p+1}\in V\mathrm{D}$ or, which is the same, differential
$k$-forms $\omega$ locally of the form
\[
\omega=\sum_{l\geq0}\omega_{i_{1}\cdots i_{p+l}a_{1}\cdots a_{k-p-l}}%
dx^{i_{1}}\cdots dx^{i_{p+l}}dy^{a_{1}}\cdots dy^{a_{k-p-l}},
\]
$\ldots,\omega_{i_{1}\cdots i_{p+l}a_{1}\cdots a_{k-p-l}},\ldots$ being local
functions on $P$, $i_{1},\ldots,i_{p+l}=1,\ldots,n$, $a_{1},\ldots
,a_{k-p-l}=1,\ldots,m$.

Denote by $V\!\Lambda(P,\alpha)=\bigoplus_{k}V\!\Lambda^{k}(P,\alpha)$ (or
simply $V\!\Lambda=\bigoplus_{k}V\!\Lambda^{k}$) the quotient differential
algebra $\Lambda(P)/\Lambda_{1}$, with $d^{V}:V\!\Lambda\longrightarrow
V\!\Lambda$ its differential and with $p^{V}:\Lambda(P)\ni\omega
\longmapsto\omega^{V}:= \omega + \Lambda_1 \in V\!\Lambda$ the projection onto the quotient. Notice
that $d^{V}$ is $A_{0}$-linear. An element $\rho^{V}$ in $V\!\Lambda^{k}$ is
locally of the form
\[
\rho^{V}=\rho_{a_{1}\cdots a_{k}}d^{V}\!y^{a_{1}}\cdots d^{V}\!y^{a_{k}},
\]
$\ldots,\rho_{a_{1}\cdots a_{k}},\dots$ being local functions on $P$, and $d^{V}\!%
\rho^{V}$ is locally given by
\[
d^{V}\!\rho^{V}=\partial_{a}\rho_{a_{1}\cdots a_{k}}d^{V}\!y^{a}d^{V}\!y^{a_{1}%
}\cdots d^{V}\!y^{a_{k}}=\partial_{\lbrack a}\rho_{a_{1}\cdots a_{k}]}d^{V}\!%
y^{a}d^{V}\!y^{a_{1}}\cdots d^{V}\!y^{a_{k}}.
\]
Clearly, $V\!\Lambda^{1}$ is the dual $A$-module of $V\mathrm{D}$ and
$V\!\Lambda$ its exterior algebra. In particular, elements in $V\!\Lambda$ may
be interpreted as multilinear, skew-symmetric forms on $V\mathrm{D}$.

Denote by $\overline{\Lambda}(P,\alpha)=\bigoplus_{k}\overline{\Lambda}{}%
^{k}(P,\alpha):=\bigoplus_{k}\Lambda_{k}^{k}\subset\Lambda(P)$ (or simply
$\overline{\Lambda}=\bigoplus_{k}\overline{\Lambda}{}^{k}$) the sub-algebra
generated by $\Lambda_{1}^{1}$. An element $\omega\in\overline{\Lambda}{}^{k}$
is locally of the form
\[
\omega=\omega_{i_{1}\cdots i_{k}}dx^{i_{1}}\cdots dx^{i_{k}}.
\]
Notice that $\overline{\Lambda}$ is naturally isomorphic to $A\otimes_{A_{0}%
}\Lambda(M)$ as an $A$-algebra.

For any $p$, the quotient (graded)
differential module $E_{0}^{p,\bullet}\equiv E_{0}^{p,\bullet}%
(P,\alpha):=\Lambda_{p}/\Lambda_{p+1}$\footnote{This last notation is
motivated by the fact that $A$--modules $E_{0}^{p,\bullet}$ are columns of the
first term of the (cohomological) Leray-Serre spectral sequence of the fiber
bundle $\alpha$ (see \cite{mvv??}).} is naturally isomorphic to $V\!\Lambda
\otimes_{A}\overline{\Lambda}{}^{p}$ (or, which is the same, $V\!\Lambda
\otimes_{A_{0}}\Lambda^{p}(M)$) via the correspondence
\begin{equation}
E_{0}^{p,q}\ni\omega+\Lambda_{p+1}^{p+q}\longmapsto\varpi\in V\!\Lambda
^{q}\otimes_{A}\overline{\Lambda}{}^{p}, \label{EqXX}
\end{equation}
well defined by putting 
\[
\varpi(Y_{1},\ldots,Y_{q}):=(i_{Y_{q}}\circ
\cdots\circ i_{Y_{1}})(\omega)\in\overline{\Lambda}{}^{p},
\]
$Y_{1}%
,\ldots,Y_{q}\in V\mathrm{D}$. In the following we denote by $E_{0}^{p,q}$ the $q$th homogeneous piece
of $E_{0}^{p,\bullet}$, $q\in\mathbb{Z}$. According to the above said, $V\!\Lambda
\otimes_{A}\overline{\Lambda}{}$ (or, which is the same, $V\!\Lambda
\otimes_{A_{0}}\Lambda(M)$) is the graded object associated with the
filtration $\Lambda(P)\supset\Lambda_{1}\supset\cdots\supset\Lambda_{p}%
\supset\cdots$. As we will see in the next subsection, a connection in
$\alpha$ allows one to identify such filtration with its graded object.

Let us now focus on the ideals $\Lambda_{n-1}$ and $\Lambda_{n}$. Put
$d^{n}x:=dx^{1}\cdots dx^{n}$ and $d^{n-1}x_{i}:=i_{\partial_{i}}d^{n}x$, so
that $dx^{j}d^{n-1}x_{i}=\delta_{i}^{j}d^{n}x$, $i,j=1,\ldots,n$. Then an
element $\omega\in\Lambda_{n-1}^{q+n-1}$ (resp.{} $\omega\in\Lambda_{n}^{q+n-1}%
$) is locally in the form
\[
\omega=\omega_{a_{1}\ldots a_{q}}^{i}dy^{a_{1}}\cdots dy^{a_{q}}d^{n-1}%
x_{i}+\omega_{a_{1}\ldots a_{q-1}}dy^{a_{1}}\cdots dy^{a_{q-1}}d^{n}x
\]
(resp.
\[
\omega=\omega_{a_{1}\ldots a_{q-1}}dy^{a_{1}}\cdots dy^{a_{q-1}}
d^{n}x),
\]$\ldots,\omega_{a_{1}\ldots a_{q}}^{i},\ldots,\omega_{a_{1}\ldots
a_{q-1}},\ldots$ being local functions on $P$. In particular $\Lambda_{n-1}^{q+n-1}$
(resp.{} $\Lambda_{n}^{q+n-1}$) is the module of sections of an $\left[
n\tbinom{q}{m}+\tbinom{q-1}{m}\right]  $(resp.{} $\tbinom{q-1}{m}$%
)-dimensional vector bundle over $P$. In few lines we will provide an
alternative description of $\Lambda_{n-1}$ and $\Lambda_{n}$ (Theorem \ref{ThIso}). In our opinion, such description
is more suitable for a better understanding of the role of $\Lambda_{n-1}$ and $\Lambda_{n}$ in
first order field theories (see, for instance, \cite{gim98}).

\subsection{Affine Forms}

Let $\nabla\in C \equiv C(P,\alpha)$. Recall, preliminarily, that $C$ is an affine space
modeled over the $A$-module $\overline{\Lambda}{}^{1}\otimes_{A}V\mathrm{D,}$
or, which is the same, $\Lambda^{1}(M)\otimes_{A_{0}}V\mathrm{D}$. The connection $\nabla$
allows one to split the tangent bundle $TP$ to $P$ into its vertical part $V\!P$
and a horizontal part $H_{\nabla}P$. denote by $H_{\nabla}\mathrm{D}%
(P,\alpha)\subset\mathrm{D}(P)$ (or simply $H_{\nabla}\mathrm{D}%
\subset\mathrm{D}(P)$) the submodule of $\nabla$-horizontal vector fields. An
element $X\in H_{\nabla}\mathrm{D}$ is locally in the form $X=X^{i}\nabla_{i}%
$, where $\nabla_{i}:=\partial_{i}+\nabla_{i}^{a}\partial_{a}$, $i=1,\ldots
,n$. Splitting 
\begin{equation}
\mathrm{D}(P)=V\mathrm{D}\oplus H_{\nabla}\mathrm{D} \label{EqSplit}
\end{equation}
determines a splitting of the de Rham differential $d:\Lambda
(P)\longrightarrow\Lambda(P)$ into a horizontal part $d_{\nabla}%
:\Lambda(P)\longrightarrow\Lambda(P)$, and a vertical part $d_{\nabla}%
^{V}:\Lambda(P)\longrightarrow\Lambda(P)$, $d=d_{\nabla}+d_{\nabla}^{V}$,
where $d_{\nabla}$ (resp.{} $d_{\nabla}^{V}$) is the Lie derivative along the
horizontal-form valued vector field (resp.{} the form valued vertical vector field)
$H_{\nabla}:A\longrightarrow\overline{\Lambda}{}^{1}(P)$ (resp.{} $V_{\nabla
}:A\longrightarrow\Lambda^{1}(P)$) determined by $\nabla$. $H_{\nabla}$ (resp.{} 
$V_{\nabla}$) is locally given by $H_{\nabla}=dx^{i}\nabla_{i}$ (resp.{} 
$V_{\nabla}=(dy^{a}-\nabla_{i}^{a}dx^{i})\partial_{a}$). Notice that
$(\Lambda(P),d_{\nabla}^{V},d_{\nabla})$ is not a bi-complex unless $\nabla$
is flat. Splitting (\ref{EqSplit}) also determines an isomorphism $\phi_{\nabla}:V\!\Lambda
\otimes_{A}\overline{\Lambda}{}\longrightarrow\Lambda(P)$ locally given by
\[
\phi_{\nabla}(d^{V}\!y^{a_{1}}\cdots d^{V}\!y^{a_{q}}\otimes dx^{i_{1}}\cdots
dx^{i_{q}})=d_{\nabla}^{V}y^{a_{1}}\cdots d_{\nabla}^{V}y^{a_{q}}dx^{i_{1}%
}\cdots dx^{i_{q}}.
\]
In particular, for any $q,p$, there is an obvious projection $\mathfrak{p}_{\nabla
}^{q,p}:\Lambda(P)\longrightarrow V\!\Lambda^{q}\otimes_{A}\overline{\Lambda
}{}{}^{p}$.

For any $k\geq0$ put 
\[
{}^\prime \Omega{}^{k+1}:=\mathrm{Aff}_{A}%
(C,V\!\Lambda^{k}\otimes_{A}\overline{\Lambda}{}^{n}).
\] An element
${}^\prime{\vartheta}\in{}^\prime \Omega{}^{k+1}$ is locally given by
\[
{}^\prime{\vartheta}(\nabla)=({}^\prime{\vartheta}{}_{a,a_{1}\cdots a_{k}}%
^{i}\nabla_{i}^{a}+{}^\prime{\vartheta}{}_{a_{1}\cdots a_{k}})d^{V}\!y^{a_{1}%
}\cdots d^{V}\!y^{a_{k}}\otimes d^{n}x,\quad\nabla\in C,
\]
$\ldots,{}^\prime{\vartheta}{}_{a,a_{1}\cdots a_{k}}^{i},\ldots,{}^\prime
{\vartheta}{}_{a_{1}\cdots a_{k}},\ldots$ local functions on $P$. The linear
part ${}^\prime\underline{{\vartheta}}$ of an element ${}^\prime{\vartheta}\in$
${}^\prime \Omega{}^{k+1}$ is an element in the $A$-module
\[
\mathrm{Hom}_{A}(\overline{\Lambda}{}^{1}\otimes_{A}V\mathrm{D},V\!\Lambda
^{k}\otimes_{A}\overline{\Lambda}{}^{n})  \simeq\mathrm{Hom}_{A}%
(V\mathrm{D},V\!\Lambda^{k}\otimes_{A}\overline{\Lambda}{}^{n-1}),
\]
where we identified $V\!\Lambda^{k}\otimes\overline{\Lambda}{}^{n-1}$ and
$\mathrm{Hom}_{A}(\overline{\Lambda}{}^{1},V\!\Lambda^{k}\otimes
\overline{\Lambda}{}^{n})$ via the isomorphism
\[
V\!\Lambda^{k}\otimes\overline{\Lambda}{}^{n-1}\ni\sigma\otimes\rho
\longmapsto\varphi_{\sigma\otimes\rho}\in\mathrm{Hom}_{A}(\overline{\Lambda}%
{}^{1},V\!\Lambda^{k}\otimes\overline{\Lambda}{}^{n}),
\]
$\sigma\in V\!\Lambda^{k}$, $\rho\in\overline{\Lambda}{}^{n-1}$, defined by
putting
\[
\varphi_{\sigma\otimes\rho}(\eta):=(-)^{k}\sigma\otimes\eta\rho\in
V\!\Lambda^{k}\otimes\overline{\Lambda}{}^{n},\quad\eta\in\overline{\Lambda}%
{}^{1}.
\]
Put $\underline{\Omega}^{k+1} \equiv \underline{\Omega}^{k+1}(P,\alpha) := V\!\Lambda^{k+1}\otimes_{A}\overline{\Lambda}{}^{n-1}$. Similarly as above, $\underline{\Omega}^{k+1}$ can be embedded into
$\mathrm{Hom}_{A}(V\mathrm{D},V\!\Lambda^{k}\otimes_{A}\overline{\Lambda}%
{}^{n-1})$ via the correspondence
\begin{equation}
\underline{\Omega}^{k+1}\ni\sigma^{\prime}%
\otimes\rho\longmapsto\varphi_{\sigma^{\prime}\otimes\rho}^{\prime}%
\in\mathrm{Hom}_{A}(V\mathrm{D},V\!\Lambda^{k}\otimes_{A}\overline{\Lambda}%
{}^{n-1}), \label{embed}
\end{equation}
$\sigma^{\prime}\in V\!\Lambda^{k}$, $\rho\in\overline{\Lambda}{}^{n-1}$,
defined by putting
\[
\varphi_{\sigma^{\prime}\otimes\rho}^{\prime}(Y):=i_{Y}\sigma^{\prime}\otimes
\rho\in V\!\Lambda^{k}\otimes_{A}\overline{\Lambda}{}^{n-1},\quad
Y\in V\mathrm{D}. 
\]
In the following we will understand embedding (\ref{embed}).\\ 
Put also $\Omega^{0}\equiv\Omega^{0}(P,\alpha):=\overline{\Lambda}{}^{n-1}$,
$\underline{\Omega}^{0}\equiv{}\underline{\Omega}^{0}(P,\alpha):=\Omega^{0}$,
and for $k\geq0$,
\[
\Omega^{k+1}\equiv\Omega^{k+1}(P,\alpha):=\{\vartheta\in{}^\prime \Omega%
{}^{k+1}\;|\;\underline{\vartheta}\in\underline{\Omega}^{k+1}\},
\]
$\Omega\equiv\Omega(P,\alpha):=\bigoplus_{q\geq0}\Omega^{q}$ and
$\underline{\Omega}\equiv\underline{\Omega}(P,\alpha):=\bigoplus_{q\geq
0}\underline{\Omega}^{q}$. Elements in $\Omega^{k}$ will be called
\emph{affine }$k$\emph{-forms} over $\alpha$, $k\geq 0$. It is easy to show that an
element $\vartheta\in{}^\prime \Omega{}^{k+1}$ is an affine $(k+1)$-form iff
it is locally given by
\[
\vartheta(\nabla)=(\vartheta_{aa_{1}\cdots a_{k}}^{i}\nabla_{i}^{a}%
+\vartheta_{a_{1}\cdots a_{k}})d^{V}\!y^{a_{1}}\cdots d^{V}\!y^{a_{k}}\otimes
d^{n}x,\quad\nabla\in C.
\]
$\ldots,\vartheta_{aa_{1}\cdots a_{k}}^{i},\ldots,\vartheta_{a_{1}\cdots
a_{k}},\ldots$ being local functions on $P$ such that $\vartheta_{aa_{1}\cdots
a_{k}}^{i}=\vartheta_{\lbrack aa_{1}\cdots a_{k}]}^{i}$, $i=1,\ldots,n$,
$a,a_{1},\ldots,a_{k}=1,\ldots,m$.

According to the above said, the linear
part $\underline{\vartheta}\in\underline{\Omega}^{k+1}$ of $\vartheta$ is
implicitly defined by the formula
\[
\vartheta(\nabla+\eta\otimes Y)(Y_{1},\ldots,Y_{k})-\vartheta(\nabla
)(Y_{1},\ldots,Y_{k})=(-)^{k}\eta\cdot\underline{\vartheta}(Y,Y_{1}%
,\ldots,Y_{k})\in\overline{\Lambda}{}^{n},
\]
$\nabla\in C$, $\eta\in\overline{\Lambda}{}^{n-1}$, $Y,Y_{1},\ldots,Y_{k}\in
V\mathrm{D}$, and it is locally given by
\begin{equation}
\underline{\vartheta}=\tfrac{(-)^{k}}{k+1}\vartheta_{a_{1}\cdots a_{k+1}}%
^{i}d^{V}\!y^{a_{1}}\cdots d^{V}\!y^{a_{k+1}}\otimes d^{n-1}x_{i}. \label{Eq4}%
\end{equation}

\subsection{Affine Forms and Differential Forms\label{SecAff}}

Let $\Omega_{0}(P,\alpha)=\bigoplus_{q\geq0}$ $\Omega_{0}^{q}(P,\alpha)$ (or,
simply, $\Omega_{0}\equiv\bigoplus_{q\geq0}$ $\Omega_{0}^{q}$) be the kernel
of the projection $\Omega\ni\vartheta\longmapsto\underline{\vartheta}%
\in\underline{\Omega}$. Clearly, $\Omega_{0}^{q}$ is canonically isomorphic to
$V^{q-1}\Lambda\otimes_{A}\overline{\Lambda}{}{}^{n}$ for $q>0$ (and in the
following we will understand such isomorphism), while $\Omega_{0}^{0}=0$.
Moreover, $\Omega^{q}$ (resp.{} $\Omega_{0}^{q}$) is the module of sections of
an $[n\tbinom{q}{m}+\tbinom{q-1}{m}]$ (resp.{} $\tbinom{q-1}{m}$)-dimensional
vector bundle over $P$.

\begin{thrm}\label{ThIso}
There are canonical isomorphisms of $A$-modules
\begin{align*}
\iota_{0,q}  &  :\Lambda_{n}^{q+n-1}\longrightarrow\Omega_{0}^{q},\\
\iota_{q}  &  :\Lambda_{n-1}^{q+n-1}\longrightarrow\Omega^{q},\\
\underline{\iota}{}_{q}  &  :E_{0}^{n-1,q}\longrightarrow\underline{\Omega
}^{q},
\end{align*}
$q\geq0$, such that diagram
\begin{equation}%
\begin{array}
[c]{c}%
\xymatrix{ 0  \ar[r] & \Lambda_n  \ar[d]_-{\iota_0}  \ar[r]  & \Lambda_{n-1}  \ar[d]_-{\iota}  \ar[r] & E_0^{n-1}  \ar[d]_-{\underline{\iota}}  \ar[r] & 0 \\
0   \ar[r] & \Omega_0                              \ar[r] & \Omega                                    \ar[r] & \underline{\Omega}                                                                   \ar[r] & 0
}
\end{array}
\label{Diag1}%
\end{equation}
commutes, where $\iota_{0}:=\bigoplus_{q}\iota_{0,q}$, $\iota:=\bigoplus
_{q}\iota_{q}$ and $\underline{\iota}:=\bigoplus_{q}\underline{\iota}{}_{q}$.
\end{thrm}

\begin{proof}
Let $q>0$. First of all, denote by $\underline{\iota}{}_{q}:E_{0}%
^{n-1,q}\longrightarrow\underline{\Omega}^{q}$ the already mentioned natural
isomorphism (\ref{EqXX}) and notice that for any $\omega\in\Lambda_{n}^{q+n-1}$ and
$Y_{1},\ldots,Y_{q-1}\in V\mathrm{D}$, $(i_{Y_{1}}\circ\cdots\circ i_{Y_{q-1}%
})(\omega)\in\overline{\Lambda}{}^{n}$. Therefore, it is well defined an
element $\iota_{0,q}(\omega)\in\Omega_{0}^{q}$ by putting $\iota_{0,q}%
(\omega)(Y_{1},\ldots,Y_{q-1}):=(i_{Y_{1}}\circ\cdots\circ i_{Y_{q-1}}%
)(\omega)\in\overline{\Lambda}{}^{n}$, $Y_{1},\ldots,Y_{q-1}\in V\mathrm{D}$.
Moreover, the correspondence $\Lambda_{n}^{q+n-1}\ni\omega\longmapsto
\iota_{0,q}(\omega)\in\Omega_{0}^{q}$ is an isomorphism of $A$-modules.
Indeed, let $\omega\in\Lambda_{n}^{q+n-1}$ and $(i_{Y_{1}}\circ\cdots\circ
i_{Y_{q-1}})(\omega)=0$ for all $Y_{1},\ldots,Y_{q-1}\in V\mathrm{D}$, then
$\omega\in\Lambda_{n+1}^{q+n-1}=\boldsymbol{0}$, so that $\iota_{0,q}$ is injective. Moreover, $\Lambda_{n}^{q+n-1}$ and
$\Omega_{0}^{q}$ are locally free $A$-modules of the same local dimension. We
conclude that
\[
\iota_{0}:=%
{\textstyle\bigoplus\nolimits_{q}}
\iota_{0,q}:\Lambda_{n}\longrightarrow\Omega_{0}%
\]
is a canonical isomorphism of $A$-modules as well, sending
$\Lambda_{n}^{q+n-1}$ into $\Omega_{0}^{q}$, $q\geq0$. Finally, if $\omega
\in\Lambda_{n}^{q+n-1}$ is locally given by
\[
\omega=\omega_{a_{1}\ldots a_{q-1}}dy^{a_{1}}\cdots dy^{a_{q-1}}d^{n}x,
\]
then $\iota_{0}(\omega)\in\Omega_{0}$ is locally given by $\iota_{0}%
(\omega)=\omega_{a_{1}\ldots a_{q-1}}d^{V}\!y^{a_{1}}\cdots d^{V}\!y^{a_{q-1}%
}\otimes d^{n}x$.

Now, for $\omega\in\Lambda_{n-1}^{q+n-1}$ and $\nabla\in C$ put
\[
\iota_{q}(\omega)(\nabla):=\mathfrak{p}_{\nabla}^{q-1,n}(\omega).
\]
If $\omega$ is locally given by
\[
\omega=\omega_{a_{1}\ldots a_{q}}^{i}dy^{a_{1}}\cdots dy^{a_{q}}d^{n-1}%
x_{i}+\omega_{a_{1}\ldots a_{q-1}}dy^{a_{1}}\cdots dy^{a_{q-1}}d^{n}x,
\]
then
\begin{eqnarray*}
\omega &  =& \omega_{a_{1}\ldots a_{q}}^{i}(d_{\nabla}^{V}+d_{\nabla})(y^{a_{1}%
})\cdots(d_{\nabla}^{V}+d_{\nabla})(y^{a_{q}})d^{n-1}x_{i}\\
& &{} +\omega_{a_{1}\ldots a_{q-1}}d_{\nabla}^{V}y^{a_{1}}\cdots d_{\nabla}%
^{V}y^{a_{q-1}}d^{n}x\\
&  =& \omega_{a_{1}\ldots a_{q-1}}d_{\nabla}^{V}y^{a_{1}}\cdots d_{\nabla}%
^{V}y^{a_{q-1}}d^{n}x + \omega_{a_{1}\ldots a_{q}}^{i}d_{\nabla}^{V}y^{a_{1}}\cdots d_{\nabla}%
^{V}y^{a_{q}}d^{n-1}x_{i} \\
& &{} +\sum_{s}(-)^{p-s}\omega_{a_{1}\ldots a_{q}}^{i}\nabla_{j}^{a_{s}}%
d_{\nabla}^{V}y^{a_{1}}\cdots\widehat{d_{\nabla}^{V}y^{a_{s}}}\cdots
d_{\nabla}^{V}y^{a_{q}}dx^{j}d^{n-1}x_{i}\\
&  =&\omega^{q,n-1}+[q(-)^{q-1}\omega_{aa_{1}\ldots a_{q-1}}^{i}\nabla_{i}%
^{a}+\omega_{a_{1}\ldots a_{q-1}}]d_{\nabla}^{V}y^{a_{1}}\cdots d_{\nabla}%
^{V}y^{a_{q-1}}d^{n}x
\end{eqnarray*}
where a cap \textquotedblleft$\widehat{\quad}$\textquotedblright\ denotes
omission of the factor below it, and $\omega^{q,n-1}\in\Lambda(P)$ is a
suitable form such that $\mathfrak{p}_{\nabla}^{q-1,n}(\omega^{q,n-1})=0$.
Therefore, locally
\begin{align}
\iota_{q}(\omega)(\nabla)  &  =\mathfrak{p}_{\nabla}^{q-1,n}(\omega
)\nonumber\\
&  =[q(-)^{q-1}\omega_{aa_{1}\ldots a_{q-1}}^{i}\nabla_{i}^{a}+\omega
_{a_{1}\ldots a_{q-1}}]d^{V}\!y^{a_{1}}\cdots d^{V}\!y^{a_{q-1}}\otimes d^{n}x.
\label{Eq5}%
\end{align}
This shows simultaneously that $\iota_{q}(\omega)$ is affine, that it is in
$\Omega^{q}$ and that $\iota_{q}$ is injective. Since $\Lambda_{n-1}^{p+n}$
and $\Omega^{q}$ are locally free $A$-modules of the same local dimension,
then the correspondence $\iota_{q}:$ $\Lambda_{n-1}^{q+n-1}\ni\omega
\longmapsto\iota_{q}(\omega)\in\Omega^{q}$ is an isomorphism. Commutativity of
diagram (\ref{Diag1}) immediately follows from local formulas (\ref{Eq4}) and
(\ref{Eq5}).
\end{proof}

Notice that isomorphism $\iota$ generalizes considerably the well known isomorphism $\Lambda_{n-1}^{n}%
\simeq\mathrm{Aff}_{A}(C,\overline{\Lambda}{}^{n})$ \cite{gim98}.

Finally, let $\pi:E\longrightarrow M$ be a fiber bundle and let $\ldots,q^{A},\ldots$ be fiber coordinates on $E$.
Notice that $\Omega^{1}(E,\pi)$ (resp.{} $\underline{\Omega}^{1}(E,\pi)$) is the
$C^{\infty}(E)$-module of sections of a vector bundle $\mu_{0}\pi
:\mathscr{M}\pi\longrightarrow E$ (resp.{} $\tau_{0}^{\dag}\pi:J^{\dag}%
\pi\longrightarrow E$). Recall that there is a distinguished element $\Theta$
in $\Omega^{1}(\mathscr{M}\pi,\mu\pi)$ (resp.{} $\underline{\Theta}\in
\underline{\Omega}^{1}(J^{\dag}\pi,\tau^{\dag}\pi)$), with $\mu\pi := \pi \circ \mu_0 \pi$ (resp.{} $\tau^\dag \pi:= \pi \circ \tau_0 ^\dag \pi$), the tautological one
\cite{gim98}, which in standard coordinates \[ \ldots,x^{i},\ldots,q^{A}%
,\ldots,p_{A}^{i},\ldots,p \] on $\mathscr{M}\pi$ (resp.{} $\ldots,x^{i}%
,\ldots,q^{A},\ldots,p_{A}^{i},\ldots$ on $J^{\dag}\alpha$) is given by
\[
\Theta=p_{A}^{i}dq^{A}d^{n-1}x_{i}-pd^{n}x\quad\text{(resp.{} }\underline
{\Theta}=p_{A}^{i}d^{V}\!q^{A}\otimes d^{n-1}x_{i}\text{).}%
\]

\section{Affine Form Calculus\label{SecCalculus}}

\subsection{Natural Operations with Affine Forms}

In this section we derive the main formulas of calculus on affine forms. Such
formulas will turn useful in generalizing proofs from the context of
Hamiltonian systems to the context of PD Hamiltonian systems (see Section \ref{SecNoether}).

Let $\alpha:P\longrightarrow M$ be as in the previous section. Isomorphism
$\iota$ (resp.{} $\iota_{0}$, $\underline{\iota}$) can be used to
\textquotedblleft transfer structures\textquotedblright\ from $\Lambda_{n-1}$
(resp.{} $\Lambda_{n}$, $E_{0}^{n-1}$) to $\Omega$ (resp.{} $\Omega_{0}$,
$\underline{\Omega}$) and back. As an instance, notice that $\Omega$ has
got a natural structure of $\Lambda(P)$-module given by
\[
\lambda\vartheta:=\iota(\lambda\omega),
\]
$\lambda\in\Lambda(P)$, $\vartheta=\iota(\omega)\in\Omega$, $\omega\in
\Lambda_{n-1}$. Moreover, $\Omega$ is generated by $\Omega^{0}$ as a
$\Lambda(P)$-module. Similarly, $\Lambda_{n}$ (resp.{} $E_{0}^{n-1}$) has a
structure of $V\!\Lambda$-module given by
\[
\lambda^{V}\omega_{0}:=\iota_{0}^{-1}(\lambda^{V}\rho^{V}\otimes\nu)\text{
(resp.{} }\lambda^{V}\underline{\omega}:=\iota_{0}^{-1}(\lambda^{V}\rho
^{V}\otimes\sigma)\text{),}%
\]
$\lambda^{V}\in V\!\Lambda$, $\omega_{0}=\iota_{0}^{-1}(\rho^{V}\otimes\nu)$
(resp.{} $\underline{\omega}=\underline{\iota}^{-1}(\rho^{V}\otimes\sigma)$),
$\rho^{V}\in V\!\Lambda$, $\nu\in\overline{\Lambda}{}^{n}$ (resp.{} $\sigma\in$
$\overline{\Lambda}{}^{n-1}$), so that $\rho^{V}\otimes\nu\in V\!\Lambda
\otimes_{A}\overline{\Lambda}{}^{n}=\Omega_{0}$ (resp.{} $\rho^{V}\otimes
\sigma\in V\!\Lambda\otimes_{A}\overline{\Lambda}{}^{n-1}=\underline{\Omega}%
$). Clearly, $\Lambda_{n}$ (resp.{} $E_{0}^{n-1}$) is generated by
$\overline{\Lambda}{}^{n}$ (resp.{} $\overline{\Lambda}{}^{n-1}$) as a
$V\!\Lambda$-module. Finally, the presented structures are compatible in the
sense that for $\omega_{0}\in\Lambda_{n}$, $\omega\in\Lambda_{n-1}$ and
$\lambda\in\Lambda(P)$, we have
\[
\lambda^{V}\omega_{0}=\lambda\omega_{0}\text{ and }\underline{\lambda\omega
}=\lambda^{V}\underline{\omega}.
\]
As a last instance of how to use isomorphisms in (\ref{Diag1}) to transfer a
structure from one space to the other we define the insertion of a connection
$\nabla\in C$ into a differential form $\omega\in\Lambda_{n}$ as follows
\[
i_{\nabla}\omega:=\iota_{0}^{-1}(\vartheta(\nabla))=(\iota_{0}^{-1}%
\circ\mathfrak{p}_{\nabla}^{q-1,n})(\omega)\in\Lambda_n,
\]
$\vartheta=\iota(\omega)\in \Omega$. Notice that the just defined insertion of a
connection in an element $\omega\in\Lambda_{n}$ has been already discussed in
\cite{emr96}. In the following we will always understand isomorphisms $\iota$,
$\iota_{0}$, $\underline{\iota}$.

Notice that $\underline{\Omega}$ inherits many operations from $\Omega$.
Indeed, let $\nabla\in C$, $Z\in\overline{\Lambda}{}^{1}\otimes_{A}%
V\mathrm{D} \subset \Lambda(P) \otimes_A \mathrm{D}(P)$, $Y\in V\mathrm{D}$, $X\in\mathrm{D}_{V}$, $q\geq0$. Then

\begin{itemize}
\item $i_{Z}(\Omega)\subset\Omega_{0}$ and $i_{Z}(\Omega_{0})=0$ so that an operator, which, abusing the notation, we again denote by
$i_{Z}:\underline{\Omega}\longrightarrow\Omega_{0}$, is well defined via the formula
\[
i_{Z}\underline{\omega}:=i_{Z}\omega\in\Omega_{0},
\]
$\omega\in\Omega$. Moreover, it is easy to show that 
\[
i_{Z}\underline{\omega
}=i_{\nabla+Z}\omega-i_{\nabla}\omega.
\]
Finally, for $Z=\eta\otimes Y_{1}$,
and $\underline{\omega}=\rho^{V}\otimes\sigma$, $\eta\in\overline{\Lambda}%
{}^{1}$, $Y_{1}\in V\mathrm{D}$, $\rho^{V}\in V\!\Lambda^{q}$ and $\sigma
\in\overline{\Lambda}{}^{n-1}$, we have
\[
i_{Z}\underline{\omega}=(-)^{q-1}i_{Y_{1}}\rho^{V}\otimes\eta\sigma.
\]

\item $i_{Y}(\Omega)\subset\Omega$ (resp.{} $L_{X}(\Omega)\subset\Omega$) and
$i_{Y}(\Omega_{0})\subset\Omega_{0}$ (resp.{} $L_{X}(\Omega_{0})\subset
\Omega_{0}$) so that the quotient map, which, abusing the
notation, we again denote by $i_{Y}:\underline{\Omega}\longrightarrow
\underline{\Omega}$ (resp.{} $L_{X}:\underline{\Omega}\longrightarrow
\underline{\Omega}$), is well defined via the formula
\[
i_{Y}\underline{\omega}:=\underline{i_{Y}\omega}\in\underline{\Omega}\text{
(resp.{} }L_{X}\underline{\omega}=\underline{L_{X}\omega}\in\underline{\Omega
}\text{)}.
\]
Finally, for $\underline{\omega}=\rho^{V}\otimes\sigma$, $\rho^{V}\in
V\!\Lambda^{q}$ and $\sigma\in\overline{\Lambda}{}^{n-1}$, we have
\[
i_{Y}\underline{\omega}=i_{Y}\rho^{V}\otimes\sigma.
\]

\item $d_{\nabla}(\Omega)\subset\Omega_{0}$ and $d_{\nabla}(\Omega_{0})=0$ so
that an operator, which, abusing the notation, we again
denote by $d_{\nabla}:\underline{\Omega}\longrightarrow\Omega_{0}$, is well defined via the
formula
\[
d_{\nabla}\underline{\omega}:=d_{\nabla}\omega\in\Omega_{0},
\]
$\omega\in\Omega$.
\end{itemize}

\begin{remark}
\label{RemPoint}Notice that the insertion $i_{\nabla}\omega$, being affine in
$\nabla$, is actually point wise, i.e., if $\nabla^{\prime}\in C$ is such that
$\nabla_{y}^{\prime}=\nabla_{y}\in C|_{y}=\alpha_{1,0}^{-1}(y)$ for some $y \in P$, then
$(i_{\nabla^{\prime}}\omega)_{y}=(i_{\nabla}\omega)_{y}$. Therefore, the insertion $i_{c}\omega_{y}$ of an element $c\in\alpha_{1,0}%
^{-1}(y)$, $y\in P$, into $\omega_{y}$ is well defined. Similar considerations apply to both
the above defined insertions $i_{Z}$ and $i_{Y}$. Finally, for all $y\in P$,
the projection $\Omega\longrightarrow\underline{\Omega}$ as well determines a
well defined linear map $\Omega|_{y}\ni\omega_{y}\longmapsto
\underline{\omega}_{y}\in\underline{\Omega}|_{y}$ whose kernel is $\Omega
_{0}|_{y}$.
\end{remark}

In the following we will denote by $\delta:\Omega\longrightarrow\Omega$ (resp.{} 
$\delta_{0}:\Omega_{0}\longrightarrow\Omega_{0}$) the restricted de Rham
differential, i.e., for $\omega\in\Omega$ (resp.{} $\omega_{0}\in\Omega_{0}$),
$\delta\omega:=d\omega\in\Omega$ (resp.{} $\delta_{0}\omega_{0}:=d\omega_{0}%
\in\Omega_{0}$) and with $\underline{\delta}:\underline{\Omega}\longrightarrow
\underline{\Omega}$ the quotient differential. Then, for $\omega_{0}=\rho
^{V}\otimes\alpha^{\ast}(\nu_{0})$ (resp.{} $\omega=\rho^{V}\otimes\alpha^{\ast
}(\sigma_{0})$), $\rho^{V}\in V\!\Lambda$, $\nu_{0}\in\Lambda^{n}(M)$ (resp.
$\sigma_{0}\in\Lambda^{n-1}(M)$), we have
\[
\delta_{0}\omega_{0}=d^{V}\!\rho^{V}\otimes\alpha^{\ast}(\nu_{0})\text{ (resp.
}\underline{\delta}\underline{\omega}=d^{V}\!\rho^{V}\otimes\alpha^{\ast}%
(\sigma_{0})\text{)}.
\]
In other words $\delta_{0}$ (resp.{} $\underline{\delta}$) is isomorphic to the
differential $d^{V}\otimes\mathrm{id}:V\!\Lambda\otimes_{A_{0}}\Lambda
^{n}(M)\longrightarrow V\!\Lambda\otimes_{A_{0}}\Lambda^{n}(M)$ (resp.
$d^{V}\otimes\mathrm{id}:V\!\Lambda\otimes_{A_{0}}\Lambda^{n-1}%
(M)\longrightarrow V\!\Lambda\otimes_{A_{0}}\Lambda^{n-1}(M)$).

All the above mentioned formulas can be proved by straightforward computations.

Now, let $\nabla$, $Y$ and $X$ be as above. denote by $[\![\cdot,\cdot]\!]$
the Fr\"{o}licher-Nijenhuis bracket in $\Lambda(P)\otimes_{A}\mathrm{D}(P)$.
It is easy to see that $[\![H_{\nabla},X]\!]\in\overline{\Lambda}{}^{1}%
\otimes_{A}V\mathrm{D}\subset\Lambda(P)\otimes_{A}\mathrm{D}(P)$. It holds the following

\begin{thrm}
Let $\omega\in\Omega$, then
\begin{gather}
\lbrack i_{\nabla},\delta]\omega:=(i_{\nabla}\circ\delta-\delta_{0}\circ
i_{\nabla})\omega=d_{\nabla}\omega\in\Omega_{0},\nonumber\\
\lbrack i_{\nabla},i_{Y}]\omega:=(i_{\nabla}\circ i_{Y}-i_{Y}\circ i_{\nabla
})\omega=0\in\Omega_{0},\label{Eq9}\\
\lbrack i_{\nabla},L_{X}]\omega:=(i_{\nabla}\circ L_{X}-L_{X}\circ i_{\nabla
})\omega=i_{[\![H_{\nabla},X]\!]}\omega\in\Omega_{0}.\nonumber
\end{gather}

\end{thrm}

\begin{proof}
First prove that $i_{\nabla}:\Omega\longrightarrow\Omega_{0}$ satisfies the
\textquotedblleft Leibniz rule\textquotedblright\
\begin{equation}
i_{\nabla}(\lambda\omega)=\lambda\cdot i_{\nabla}\omega+i_{H_{\nabla}}%
\lambda\cdot\omega, \label{Eq6}%
\end{equation}
$\lambda\in\Lambda(P)$, $\omega\in\Omega$. For $\rho\in\Lambda(P)$, denote
$\rho_{\nabla}^{\bullet,p}:=\sum_{q}\mathfrak{p}_{\nabla}^{q,p}(\rho)$, so
that $\rho=\sum_{p}\rho_{\nabla}^{\bullet,p}$. Notice that for $\omega
\in\Omega$ and $\lambda\in\Lambda(P)$, we have $\omega=\omega_{\nabla
}^{\bullet,n}+\omega_{\nabla}^{\bullet,n-1}$ so that
\[
i_{\nabla}(\lambda\omega)=\mathfrak{p}_{\nabla}^{\bullet,n}(\lambda
\omega)=\lambda_{\nabla}^{l,0}\omega_{\nabla}^{\bullet,n}+\lambda_{\nabla
}^{\bullet,1}\omega_{\nabla}^{\bullet,n-1}=\lambda\cdot i_{\nabla}%
\omega+\lambda_{\nabla}^{\bullet,1}\cdot\omega_{\nabla}^{\bullet,n-1}.
\]
Moreover, $i_{H_{\nabla}}\lambda=%
{\textstyle\sum\nolimits_{p}}
i_{H_{\nabla}}\lambda_{\nabla}^{\bullet,p}=%
{\textstyle\sum\nolimits_{p}}
p\lambda_{\nabla}^{\bullet,p}$, which in turn implies $\lambda_{\nabla
}^{\bullet,p}=i_{H_{\nabla}}\lambda-%
{\textstyle\sum\nolimits_{p>1}}p\lambda_{\nabla}^{\bullet,p}$. Therefore
\begin{align*}
i_{\nabla}(\lambda\omega) &  =\lambda\cdot i_{\nabla}\omega+\lambda_{\nabla}^{\bullet,1}\cdot
\omega_{\nabla}^{\bullet,n-1}\\
&  =\lambda\cdot i_{\nabla}\omega+i_{H_{\nabla}}\lambda\cdot\omega_{\nabla
}^{\bullet,n-1} - {\textstyle\sum\nolimits_{p>1}}p\lambda_{\nabla}^{\bullet,p}\omega_{\nabla}^{\bullet,n-1}\\
&  =\lambda\cdot i_{\nabla}\omega+i_{H_{\nabla}}\lambda\cdot\omega.
\end{align*}
In view of (\ref{Eq6}), the above defined operators $[i_{\nabla},\delta],$
$[i_{\nabla},i_{Y}],$ $[i_{\nabla},L_{X}]:\Omega\longrightarrow\Omega_{0}$,
satisfy analogous \textquotedblleft Leibniz rules\textquotedblright\:
\begin{align}
\lbrack i_{\nabla},\delta](\lambda\omega) & =d_{\nabla}\lambda\cdot\omega
+(-)^{l}\lambda\cdot\lbrack i_{\nabla},\delta](\omega), \nonumber \\
\lbrack i_{\nabla},i_{Y}](\lambda\omega) & =\lambda\cdot\lbrack i_{\nabla}%
,i_{Y}] (\omega),\label{EqYY} \\
\lbrack i_{\nabla},L_{X}](\lambda\omega) & =i_{[\![H_{\nabla},X]\!]}\lambda
\cdot\omega+\lambda\cdot\lbrack i_{\nabla},L_{X}](\omega). \nonumber
\end{align}

Since $\Omega$ is generated by $\overline{\Lambda}{}^{n-1}$ as a $\Lambda
(P)$-module, in view of (\ref{EqYY}), it is enough to prove (\ref{Eq9}) for $\omega
\in\overline{\Lambda}{}^{n-1}$. In this case
\begin{gather*}
\lbrack i_{\nabla},\delta]\omega=i_{\nabla}d\omega=(d\omega)_{\nabla}%
^{\bullet,n}=d_{\nabla}\omega,\\
\lbrack i_{\nabla},i_{Y}]\omega=0,\\
\lbrack i_{\nabla},L_{X}]\omega=i_{\nabla}L_{X}\omega=(L_{X}\omega)_{\nabla
}^{\bullet,n}=0=i_{[\![H_{\nabla},X]\!]}\omega.
\end{gather*}

 \end{proof}

We now discuss the interaction between affine forms and bundle morphisms. Let
$\alpha^{\prime}:P^{\prime}\longrightarrow M$ be another fiber bundle and let
$G:P\longrightarrow P^{\prime}$ be a bundle morphism. Clearly, $G$ preserves the
ideals $\Lambda_{p}$, $p\geq0$, i.e., $G^{\ast}(\Lambda_{p}(P^{\prime}%
,\alpha^{\prime}))\subset\Lambda_{p}(P,\alpha)$. In particular,
\[
G^{\ast}(\Omega(P^{\prime},\alpha^{\prime}))\subset\Omega(P,\alpha)\text{ and
}G^{\ast}(\Omega_{0}(P^{\prime},\alpha^{\prime}))\subset\Omega_{0}%
(P,\alpha)\text{.}%
\]
We conclude that the quotient map which, abusing the
notation, we again denote by $G^{\ast}:\underline{\Omega}(P^{\prime}%
,\alpha^{\prime})\longrightarrow\underline{\Omega}(P,\alpha)$, is well defined. Now, consider
$G$-compatible connections $\nabla\in C(P,\alpha)$ and $\nabla^{\prime}\in
C(P^{\prime},\alpha^{\prime})$. It is easy to show that
\begin{equation}
G^{\ast}\circ i_{\nabla^{\prime}}=i_{\nabla}\circ G^{\ast}:\Omega(P^{\prime
},\alpha^{\prime})\longrightarrow\Omega_{0}(P,\alpha). \label{Eq33}%
\end{equation}

\subsection{Cohomology}

\begin{remark}
(see \cite{mvv??}) In the following we denote by $\mathcal{F}$ the abstract
fiber of $\alpha$. Notice that, for any $q\geq0$, $V\!H^{q}\equiv
V\!H^{q}(P,\alpha):=H^{q}(V\!\Lambda,d^{V})$ is the $A_{0}$-module of sections
of a (pro-finite) vector bundle $\alpha^{q}:P^{q}\longrightarrow M$ over $M$
whose abstract fiber is $H^{q}(\mathcal{F})$. Moreover, $\alpha^{q}$ is
endowed with a canonical flat connection $\nabla^{q}$ ($\nabla^{q}$ is a
smooth analogue of Gauss-Manin connection in algebraic geometry).
Correspondingly, there is a de Rham like complex
\[
\xymatrix@C=40pt{\cdots \ar[r] & \Lambda^{p-1}\otimes_{A_{0}}V\!H^{q} \ar[r]^-{d_{1}^{p-1,q}} & \Lambda^{p}\otimes_{A_{0}}V\!H^{q} \ar[r]^-{d_{1}^{p,q}} & \cdots},
\]
whose cohomology we denote by $E_{2}^{\bullet,q}:=\bigoplus_{p}E_{2}^{p,q}$,
$E_{2}^{p,q}:=H^{p}(\Lambda(M)\otimes_{A_{0}}V\!\Lambda^{q},d_{1}^{\bullet
,q})$\footnote{Similarly as above, this last notations are motivated by the
fact that the differentials $d_{1}^{\bullet,q}$ (resp.{} the vector spaces
$E_{2}^{\bullet,q}$) are the ones in the first term (resp.{} are rows of the
second term) of the (cohomological) Leray-Serre spectral sequence of the fiber
bundle $\alpha$ \cite{mvv??}.{} }, $q\geq0$. It can be proved that, if $\alpha$
is trivial or $M$ is simply connected, then there is a (generically
non-canonical), isomorphism
\[
E_{2}^{p,q}\approx H^{p}(M)\otimes H^{q}(\mathcal{F}),\quad p,q\geq0.
\]
Finally, notice also that, for any $q\geq0$,
\begin{align*}
H^{q}(\Omega_{0},\delta_{0})  &  \simeq\Lambda^{n}(M)\otimes_{A_{0}}%
V\!H^{q},\\
H^{q}(\underline{\Omega},\underline{\delta})  &  \simeq\Lambda^{n-1}%
(M)\otimes_{A_{0}}V\!H^{q}.
\end{align*}

\end{remark}

\begin{prop}
\label{Prop1}Let $\alpha:P\longrightarrow M$ be a fiber bundle. Then, for any
$q\geq0$, there exists a short exact sequence of vector spaces
\[
\xymatrix{0 \ar[r] & \operatorname{coker} d_1^{n,q-1} \ar[r] & H^q(\Omega,\delta) \ar[r] & \ker d_1^{n-1,q} \ar[r] & 0}.
\]
In particular, $H^{q}(\Omega,\delta)\approx\operatorname{coker}d_{1}%
^{n,q-1}\oplus\ker d_{1}^{n-1,q}=E_{2}^{n,q-1}\oplus\ker d_{1}^{n-1,q}$.
\end{prop}

\begin{proof}
Consider the short exact sequence of complexes
\[
\xymatrix{0 \ar[r] & \Omega_{0} \ar[r] & \Omega \ar[r] & \underline
{\Omega} \ar[r] & 0},
\]
and the associated long sequence in cohomology
\begin{equation}
\xymatrix@C=23pt{ \cdots \ar[r] & H^{q-1}(\underline{\Omega},\underline{\delta}) \ar[r]^-{\partial}
& H^{q}(\Omega_{0},\delta_{0})                   \ar[r]
& H^{q}(\Omega,\delta)                           \ar[r]
& H^{q}(\underline{\Omega},\underline{\delta})   \ar[r]^-{\partial}
& \cdots}. \label{Eq28}%
\end{equation}
We already commented, in the above remark, that, for any $q$, $H^{q}%
(\Omega_{0},\delta_{0})$ identifies with $\Lambda^{n}(M)\otimes_{A_{0}%
}V\!H^{q}$ and $H^{q}(\underline{\Omega},\underline{\delta})$ identifies with
$\Lambda^{n-1}(M)\otimes_{A_{0}}V\!H^{q}$. Similarly, it is easy to show that
the connecting operator
\[
\partial:H^{q-1}(\underline{\Omega},\underline{\delta})\longrightarrow
H^{q}(\Omega_{0},\delta_{0})
\]
identifies with the de Rham-like differential
\[
d_{1}^{n-1,q}:\Lambda^{n-1}(M)\otimes_{A_{0}}V\!H^{q}\longrightarrow
\Lambda^{n}(M)\otimes_{A_{0}}V\!H^{q}.
\]
The thesis then follows from exactness of (\ref{Eq28}).
 \end{proof}

\begin{cor}
\label{CorH0}If $\mathcal{F}$ is connected, then $H^{0}(\Omega,\delta
)\simeq\ker d_{M}^{n-1}$,
\[
d_{M}^{n-1}:\Lambda^{n-1}(M)\longrightarrow\Lambda^{n}(M)
\]
being the last de Rham differential of $M$.
\end{cor}

\begin{proof}
If $\mathcal{F}$ is connected $V\!H^{0}\simeq A_{0}$ and $d_{1}^{n-1,0}$
identifies with $d_{M}^{n-1}$.
 \end{proof}

\begin{cor}
\label{CorPoincare}Let $q\geq0$ and $\omega\in\Omega^{q}$ be $\delta$-closed,
i.e., $\delta\omega=0$. Then, 1) if $q=0$, $\omega$ is locally of the form
$\alpha^{\ast}(\eta)$ for some $\eta\in\Lambda^{n-1}(M)$, 2) if $q>0$, then
$\omega$ is locally $\delta$-exact, i.e., $\omega$ is locally of the form
$\delta\theta$, $\theta$ being a local element in $\Omega^{q-1}$.
\end{cor}

\begin{proof}
If $\mathcal{F}$ is contractible, then $V\!H^{q}=0$, and therefore
$H^{q}(\Omega,\delta)=0$, for all $q>0$.
 \end{proof}

Let $\omega\in\Omega$ and $\theta\in\Omega$ be such that $\omega=\delta\theta
$. Then $\theta$ will be called a \emph{potential} of $\omega$.

\section{PD Hamiltonian Systems\label{SecPDHamSys}}

\subsection{PD Hamiltonian Systems and PD Hamilton Equations}

In this section we introduce what we think should be understood as the partial
differential, i.e., field theoretic analogue of a Hamiltonian (mechanical) system on an abstract symplectic manifold.

Let $\alpha:P\longrightarrow M$ be as in the previous section and let $\omega
\in\Omega^{2}(P,\alpha)$ be such that $\delta\omega=0$. Put
\begin{align*}
\ker\omega &  :=\{Y\in V\mathrm{D}\;|\;i_{Y}\omega=0\},\quad\ker
\underline{\omega}:=\{Y\in V\mathrm{D}\;|\;i_{Y}\underline{\omega}{}=0\},\\
\operatorname{Ker}\omega &  :=\{\nabla\in C\;|\;i_{\nabla}\omega
=0\},\quad\operatorname{Ker}\underline{\omega}:=\{Z\in V\mathrm{D}\otimes
_{A}\overline{\Lambda}{}^{1}\;|\;i_{Z}\underline{\omega}{}=0\}.
\end{align*}
Since $\omega$ is closed, both $\ker\omega$ and $\ker\underline{\omega}$ are
modules of smooth sections of involutive $\alpha$-vertical distributions
$D^{\omega}$ and $\underline{D}^{\omega}$ on $P$, where, for $y\in P$,
\[
D_{y}^{\omega}:=\{\xi\in V_{y}P\;|\;i_{\xi}\omega_{y}=0\},\quad\underline
{D}_{y}^{\omega}:=\{\xi\in V_{y}P\;|\;i_{\xi}\underline{\omega}_{y}=0\}.
\]
Similarly, $\operatorname{Ker}\underline{\omega}$ is a sub-module in
$V\mathrm{D}\otimes_{A}\overline{\Lambda}{}^{1}$. As a minimal regularity
requirement, assume that $\underline{D}^{\omega}$ has got constant rank
$\underline{r}$. Then, it is easy to check that, as a consequence,
$\operatorname{Ker}\underline{\omega}$ is the module of sections of a smooth
vector bundle $\varpi:W\longrightarrow P$. For $y\in P$, denote $r(y)=\dim
D_{y}^{\omega}$. In general, $r(y)$ will change from point to point $y\in P$.
However, we are proving in brief that $r(y)$ cannot change that much. First of
all, since, obviously, $D^{\omega}\subset\underline{D}^{\omega}$, then
$r(y)\leq\underline{r}$ for all $y\in P$. Now, for $y\in P$, denote
\[
\operatorname{Ker}\omega_{y}:=\{c\in\alpha_{1,0}^{-1}(y)\;|\;i_{c}\omega
_{y}=0\}.
\]
Then, $\operatorname{Ker}\omega_{y}$ is either empty or an affine space
modeled over $\varpi^{-1}(y)$. It holds the

\begin{prop}
For any $y\in P$, $\underline{r}-r(y)\leq1$ (see also Theorem 4 of \cite{fg08}).
\end{prop}

\begin{proof}
Let $y\in P$ and suppose $r(y)<\underline{r}$. If $\xi\in\underline{D}%
_{y}^{\omega}$ then (see Remark \ref{RemPoint}) $\underline{i_{\xi}\omega_{y}%
}=i_{\xi}\underline{\omega}_{y}=0$ so that $i_{\xi}\omega_{y}\in\Omega_{0}%
^{1}|_{y}=\overline{\Lambda}{}^{n}|_{y}$. Then consider the map $\gamma
_{y}:\underline{D}_{y}^{\omega}\ni\xi\longmapsto\gamma_{y}(\xi):=i_{\xi}%
\omega_{y}\in\overline{\Lambda}{}^{n}|_{y}$. Since $r(y)<\underline{r}$,
$\gamma_{y}$ is surjective and the sequence of vector spaces $0\longrightarrow
D_{y}^{\omega}\longrightarrow\underline{D}_{y}^{\omega}\overset{\gamma_{y}%
}{\longrightarrow}\overline{\Lambda}{}^{n}|_{y}\longrightarrow0$ is exact.
Since $\overline{\Lambda}{}^{n}|_{y}$ is $1$-dimensional, it follows that
$\underline{r}-r(y)=1$.
 \end{proof}

The following proposition characterizes the case $r(y)=\underline{r}$.

\begin{prop}
\label{PropKer}Let $\omega$ be as above. Then $r(y)=\underline{r}$ iff
$\operatorname{Ker}\omega_{y}\neq\varnothing$.
\end{prop}

\begin{proof}
The result is nothing more than an application of the Rouch\'{e}-Capelli
theorem. We here propose a dual proof. Let $\xi\in V_{y}P$ be given by $\xi$
$=\xi^{a}{\partial}_a| _{y}$. Then
$\xi\in\underline{D}_{y}^{\omega}$ iff
\begin{equation}
\omega_{ab}^{i}(y)\xi^{a}=0,\quad a=1,\ldots,m,\;i=1,\ldots,n. \label{Eq12}%
\end{equation}
Similarly, $\xi\in D_{y}^{\omega}$ iff they are satisfied both (\ref{Eq12})
and
\begin{equation}
\omega_{a}(y)\xi^{a}=0. \label{Eq13}%
\end{equation}
Therefore, $D_{y}^{\omega}=\underline{D}_{y}^{\omega}$ iff Equation
(\ref{Eq13}) linearly depends on Equations (\ref{Eq12}), i.e., iff there are
real numbers $h_{i}^{b}$, $b=1,\ldots,m,i=1,\ldots,n$, such that
\[
\omega_{a}(y)=\omega_{ab}^{i}(y)h_{i}^{b},
\]
$a=1,\ldots,m$. Now, let $c\in\alpha_{1,0}^{-1}(y)$ be given by $y_{i}^{a}(c)=-\tfrac{1}%
{2}h_{i}^{a}$. Then $i_{c}\omega_{y}$ is given by
\begin{align*}
i_{c}\omega_{y}  &  =(-2\omega_{ba}^{i}(y)y_{i}^{a}(c)+\omega_{a}%
(y))dy^{a}d^{n}x\,|_{y}\\
&  =-(\omega_{ab}^{i}(y)h_{i}^{b}-\omega_{a}(y))dy^{a}d^{n}x\,|_{y}\\
&  =0.
\end{align*}

 \end{proof}

\begin{definition}
A\emph{ PD prehamiltonian system} on the fiber bundle $\alpha:P\longrightarrow
M$ is a $\delta$-closed element $\omega\in\Omega^{2}(P,\alpha)$. A
\emph{PD Hamiltonian system on }$\alpha$ is a PD prehamiltonian system
$\omega$ such that $\ker\underline{\omega}=0$ (and, therefore, $\ker\omega=0$
as well).
\end{definition}

Let $\theta\in\Omega^{1}$ be locally given by $\theta=\theta_{a}^{i}%
dy^{a}d^{n-1}x_{i}-Hd^{n}x$, $\ldots,\theta_{a}^{i},\ldots,H$ being local
functions on $P$. Then $\delta\theta$ is locally given by
\[
\delta\theta=\partial_{\lbrack a}\theta_{b]}^{i}dy^{a}dy^{b}d^{n-1}%
x_{i}-(\partial_{a}H + \partial_{i}\theta_{a}^{i})dy^{a}d^{n}x.
\]
Similarly, let $\omega\in\Omega^{2}$ and $Y\in V\mathrm{D}$ be locally given
by $\omega=\omega_{ab}^{i}dy^{a}dy^{b}d^{n-1}x_{i}+\omega_{a}dy^{a}d^{n}x$ and
$Y=Y^{a}\partial_{a}$, respectively. Then $\delta\omega$, $i_{Y}\omega$ and
$i_{Y}\underline{\omega}$ are locally given by
\begin{align*}
\delta\omega &  =\partial_{\lbrack a}\omega_{bc]}^{i}dy^{a}dy^{b}dy^{c}%
d^{n-1}x_{i}+(\partial_{i}\omega_{ab}^{i}+\partial_{\lbrack a}\omega
_{b]})dy^{a}dy^{b}d^{n}x,\\
i_{Y}\omega &  =2\omega_{ab}^{i}Y^{a}dy^{b}d^{n-1}x_{i}+\omega_{a}Y^{a}%
d^{n}x,\\
i_{Y}\underline{\omega}  &  =2\omega_{ab}^{i}Y^{a}d^{V}\!y^{b}\otimes
d^{n-1}x_{i},
\end{align*}
so that $\omega$ is a PD prehamiltonian system iff
\begin{equation}
\partial_{\lbrack a}\omega_{bc]}^{i}=0,\quad\partial_{i}\omega_{ab}%
^{i}+\partial_{\lbrack a}\omega_{b]}=0, \label{Eq31}%
\end{equation}
or, which is the same (see Corollary \ref{CorPoincare}),
\begin{equation}
\omega_{ab}^{i}=\partial_{\lbrack a}\theta_{b]}^{i},\quad\omega_{a}%
=-\partial_{a}H-\partial_{i}\theta_{a}^{i}, \label{Eq32}%
\end{equation}
for some $\ldots,\theta_{a}^{i},\ldots,H$ local functions on $P$. Moreover,
$\omega$ is a PD Hamiltonian system iff
\begin{equation}
\omega_{ab}^{i}Y^{a}=0\Longrightarrow Y^{a}=0. \label{Eq30}%
\end{equation}
In its turn (\ref{Eq30}) implies $\omega_{a}=\omega_{ab}^{i}f_{i}^{b}$ for
some $\ldots,f_{i}^{b},\ldots$ local functions on $P$ (see the proof of
Proposition \ref{PropKer}).

Let $\omega$ be a PD prehamiltonian system on $\alpha$, and let $\sigma
:U\longrightarrow P$ be a local section of $\alpha$, $U\subset M$ being an open subset.
The first jet prolongation $\dot{\sigma}:U\longrightarrow J^{1}\alpha$ of
$\sigma$ may be interpreted as a \textquotedblleft connection in $\alpha$
along $\sigma$\textquotedblright, i.e., a section of the restricted
bundle\ $\alpha_{1,0}|_{\sigma}:J^{1}\alpha|_{\sigma}\longrightarrow M$. Moreover,
elements in $\Omega|_{\sigma}$ may be interpreted as affine maps from
$C|_{\sigma}$ to $\Omega_{0}|_{\sigma}\simeq V\!\Lambda|_{\sigma}%
\otimes_{A_{0}}\Lambda^{n}(M)$ whose linear part is in $\underline{\Omega
}|_{\sigma}\simeq V\!\Lambda|_{\sigma}\otimes_{A_{0}}\Lambda^{n-1}(M)$.
Namely, an element $\lozenge\in C|_{\sigma}$ can \textquotedblleft be
inserted\textquotedblright\ into an element $\rho|_{\sigma}\in\Omega|_{\sigma
}$, $\rho\in\Omega$, giving an element $i_{\lozenge}\rho|_{\sigma}\in
\Omega_{0}|_{\sigma}$. Thus, we can search for local sections $\sigma$ of
$\alpha$ such that
\begin{equation}
i_{\dot{\sigma}}\omega|_{\sigma}=0 \label{Eq29}%
\end{equation}

\begin{definition}
Equations (\ref{Eq29}) are called the \emph{PD Hamilton equations} (of the
PD prehamiltonian system $\omega$).
\end{definition}

If $\omega$ is locally given by $\omega=\omega_{ab}^{i}dy^{a}dy^{b}%
d^{n-1}x_{i}+\omega_{a}dy^{a}d^{n}x$, then the associated PD Hamilton equations
are locally given by
\begin{equation}
2\omega_{ab}^{i}\partial_{i}y^{a}-\omega_{b}=0. \label{Eq11}%
\end{equation}
Conversely, a system of PDEs in the form (\ref{Eq11}) is a PD Hamilton
equation for some PD prehamiltonian (resp.{} PD Hamiltonian) system iff
coefficients $\ldots,\omega_{ab}^{i},\ldots,\omega_{b},\ldots$ satisfy
(\ref{Eq31}) (or, which is the same, (\ref{Eq32})) (resp.{} (\ref{Eq31}) and
(\ref{Eq30})). Notice that, in view of (\ref{Eq11}), a general PD prehamiltonian
system $\omega$ encode \textquotedblleft kinematical
information\textquotedblright, which can be identified with $\underline
{\omega}$, and \textquotedblleft dynamical information\textquotedblright,
which can be identified with the specific choice of $\omega$ in the class of
those PD Hamiltonian systems with linear part $\underline{\omega}$ (see the
comment at the end of Section \ref{SecAff}, Remark \ref{Rem1} and Example
\ref{Rem2}).

Searching for solutions of PD Hamilton equations of a PD prehamiltonian system
$\omega$, we could proceed in two steps:

\begin{enumerate}
\item \label{Prob1} search for a connection $\nabla\in\operatorname{Ker}%
\omega$,

\item search for $n$-dimensional integral submanifolds of the horizontal
distribution $H_{\nabla}P$.
\end{enumerate}

However, a solution to the first step of the above mentioned procedure exists
iff $\ker\omega=\ker\underline{\omega}$ which is not always the case.
Therefore, in general, we are led to weaken \ref{Prob1} and search for
connections $\nabla^{\prime}$ in a subbundle $P^{\prime}\subset P$ such that
$i_{\nabla^{\prime}}\omega|_{P^{\prime}}=0$. As showed in the next
proposition, there is always an \textquotedblleft algorithmic\textquotedblright\ way to find a maximal subbundle
$\breve{\alpha}:\breve{P}\longrightarrow M$ of $\alpha$ such that the affine
equation $i_{\breve{\nabla}}\omega|_{\breve{P}}=0$, $\breve{\nabla}\in
C(\breve{P},\breve{\alpha})$ admits at least one solution. We will refer to
the above mentioned \textquotedblleft algorithm\textquotedblright\ as the PD constraint algorithm (see also
\cite{gnh78,dmm96,dmm96b,d...02,d...05}).

\begin{prop}
Let $\omega$ be as above and $\operatorname{Ker}\omega=\varnothing$ (i.e.,
$D_{y}^{\omega}\neq\underline{D}_{y}^{\omega}$ for some $y\in P$). Under
suitable regularity conditions on $\omega$ (to be specified in the proof),
there exists a (maximal) subbundle $\breve{P}\subset P$ such that
$i_{\breve{\nabla}}\omega|_{\breve{P}}=0$ for some $\breve{\nabla}\in
C(\breve{P},\breve{\alpha})$.
\end{prop}

\begin{proof}
For $s=1,2,\ldots$ define recursively
\begin{align*}
P_{(s)}& :=\{y\in P_{(s-1)}\;|\;\operatorname{Ker}\omega_{y}\cap(\alpha_{(s-1)})_1
^{-1}(y)\neq\varnothing\}\subset P,\\
\alpha_{(s)} & :=\alpha|_{P_{(s)}}%
:P_{s}\longrightarrow M,
\end{align*}
where $P_{(0)}:=P$, $\alpha_{(0)}:=\alpha$ (in particular
$P_{1}=\{y\in P\;|\;r(y)=\underline{r}\}$). We assume that $\alpha_{(s)}%
:P_{(s)}\longrightarrow M$ is a smooth (closed) subbundle for all $s$
(regularity conditions). Then, for dimensional reasons, there exists
$\overline{s}$ such that $P_{(s)}=P_{(\overline{s})}$ for all $s\geq\overline{s}$.
Put $\breve{P}:=P_{(\overline{s})}$.
 \end{proof}

The subbundle $\breve{\alpha}:=\alpha|_{\breve{P}}:\breve{P}\longrightarrow M$ will be
called the \emph{constraint subbundle}. Notice that $\breve{P}$ can be empty
(for instance when $r(y)=\underline{r}-1$ for all $y\in P$) and, in this case,
 PD Hamilton equations do not possess solutions.

\begin{cor}
Let $\omega$ be a PD prehamiltonian system on $\alpha$ and let $\sigma$ be a solution
of the PD Hamilton equations. Then $\operatorname{im}\sigma\subset\breve{P}$.
\end{cor}

\begin{proof}
By induction on $s$, $\operatorname{im}\sigma\subset P_{(s)}$ for all
$s=1,2,\ldots$.
 \end{proof}

The converse of the above corollary is, a priori, only true for $n=1$.
Namely, we may wonder if for any $y\in$ $\breve{P}$ there is a solution
$\sigma$ of PD Hamilton equations such that $y\in\operatorname{im}\sigma$.
We know that there is a connection $\breve{\nabla}$ in $\breve{P}$ which is
\textquotedblleft a solution of PD Hamilton equations up to first
order\textquotedblright, i.e., $i_{\breve{\nabla}}\breve{\omega}|_{\breve{P}%
}=0$. $n$-dimensional integral manifolds of the horizontal distribution
$H_{\breve{\nabla}}\breve{P}$ determined on $\breve{P}$ by $\breve{\nabla}$
are clearly images of solutions of PD Hamilton equations. If $n=1$,
$\breve{\nabla}$ is trivially flat and Frobenius theorem guarantees that
for any $y\in\breve{P}$ there is a solution \textquotedblleft through
$y$\textquotedblright. The same is a priori untrue for $n>2$. Integrability
conditions on $H_{\breve{\nabla}}\breve{P}$ will be discussed elsewhere.

\subsection{PD Hamiltonian Systems and Multisymplectic Geometry \textit{\`{a}
l\`{a} }Forger}

Forger and Gomes have recently proposed a definition of multipresymplectic
structure on a fiber bundle \cite{fg08}. Their work aims to define such a
structure so that 1) the differential $d\Theta$ of the tautological $n$-form
$\Theta$ on the affine adjoint bundle of the first jet bundle (see the end of Section \ref{SecAff})is multisymplectic
2) every multipresymplectic structure is locally isomorphic to the pull-back
of $\Theta$ along a fibration (Darboux lemma). Since, in our opinion, I)
\cite{fg08} is the best motivated and established work about fundamentals of
multisymplectic geometry, II) abstract fiber bundles play in \cite{fg08} a
similar role as in this paper, we analyze in this subsection the relationship
between PD prehamiltonian systems and multipresymplectic structures
\textit{\`{a} l\`{a} }Forger, referring to \cite{fg08} for the main
definition. Here we just mention two of the main results of \cite{fg08} (which
can eventually be understood as definitions of polypresymplectic structure and
multipresymplectic structure on a fiber bundle, respectively)

\begin{thrm}
[Forger and Gomes I]\label{TeorForger1}Let $\alpha:P\longrightarrow M$ be a
fiber bundle, $\ldots,x^{i},\ldots$ be local coordinates on $M$, $i=1,\ldots
,n=\dim M$, and $\underline{\omega}\in\underline{\Omega}^{2}$. The form $\underline
{\omega}$ is a polypresymplectic structure on $\alpha$ iff, around every point
of $P$, there are local fiber coordinates $\ldots,q^{A},\ldots,p_{A}%
^{i},\ldots,z^{1}\ldots,z^{s}$, $A=1,\ldots,m$, $i=1,\ldots,n$ (so that $\dim
P=n+m+mn+s$) such that $\underline{\omega}$ is locally given by
\[
\omega=d^{V}\!p_{A}^{i}d^{V}\!q^{A}\otimes d^{n-1}x_{i}.
\]

\end{thrm}

\begin{thrm}
[Forger and Gomes II]\label{TeorForger2}Let $\alpha:P\longrightarrow M$ be a
fiber bundle, $\ldots,x^{i},\ldots$ be local coordinates on $M$, $i=1,\ldots
,n=\dim M$, and $\omega\in\Omega^{2}$. The form $\omega$ is a multipresymplectic
structure on $\alpha$ iff, around every point of $P$, there are local fiber
coordinates $\ldots,q^{A},\ldots,p_{A}^{i},\ldots,p,z^{1}\ldots,z^{r}$,
$A=1,\ldots,m$, $i=1,\ldots,n$ (so that $\dim P=(n+1)(m+1)+r$) such that
$\omega$ is locally given by
\[
\omega=dp_{A}^{i}dq^{A}d^{n-1}x_{i}-dpd^{n}x.
\]

\end{thrm}

\begin{prop}
\label{Prop2}Let $\omega$ be a PD prehamiltonian system on $\alpha$. The
following two conditions are equivalent.

\begin{enumerate}
\item $\underline{\omega}$ is a polypresymplectic structure and
$r(y)=\underline{r}-1$ for all $y\in P$.

\item $\omega$ is a multipresymplectic structure.
\end{enumerate}
\end{prop}

\begin{proof}
Recall that, in view of Proposition \ref{Prop1}, $\omega$ is locally $\delta
$-exact. Suppose $\underline{\omega}$ is polypresymplectic and
$r(y)=\underline{r}-1$ for all $y\in P$. Then $\underline{r}>0$ and, in view
of Theorem \ref{TeorForger1}, (around every point in $P$) there are $\alpha
$-adapted local coordinates%
\[
\ldots,x^{i},\ldots,q^{A},\ldots,p_{A}^{i},z^{0},z^{1},\ldots,z^{\underline
{r}-1},\quad
\]
such that, locally, $\underline{\omega}=d^{V}\!p_{A}^{i}d^{V}\!q^{A}\otimes
d^{n-1}x_{i}$ (in particular, $\ker\underline{\omega}$ is locally spanned by
$\ldots,{\partial}/{\partial z^{\alpha}},\ldots$). Therefore,
$\underline{\omega}=\underline{\delta}\underline{\theta}{}_{0}$, where
$\underline{\theta}{}_{0}:=p_{A}^{i}d^{V}\!q^{A}\otimes d^{n-1}x_{i}$ is a local
element of $\underline{\Omega}^{1}$. A general (local) potential of $\omega$
is then $\theta^{\prime}\in\Omega^{1}$ such that $\underline{\theta}^{\prime
}=\underline{\theta}{}_{0}+d^{V}\!\nu$, $\nu$ being a local element in
$\underline{\Omega}^{0}=\overline{\Lambda}{}^{n-1}$. The (local) potential
$\theta:=\theta^{\prime}-\delta\nu$ is locally in the form $\theta=p_{A}%
^{i}dq^{A}d^{n-1}x_{i}-pd^{n}x$, where $p$ is a local function on $P$.
Therefore $\omega$ is locally given by
\[
\omega=dp_{A}^{i}dq^{A}d^{n-1}x_{i}-dpd^{n}x.
\]
The module $\ker\omega$ is locally spanned by those local elements $Y^{\alpha}%
\tfrac{\partial}{\partial z^{\alpha}}$ in $\ker\underline{\omega}$ such that
$Y^{\alpha}\tfrac{\partial h}{\partial z^{\alpha}}=0$. Since $\ker\omega
\neq\ker\underline{\omega}$, then $\tfrac{\partial p}{\partial z^{\alpha}%
}dz^{\alpha}\neq0$. Let, for instance, be $\tfrac{\partial p}{\partial z^{0}%
}\neq0$. Then $\ldots,x^{i},\ldots,q^{A},\ldots,p_{A}^{i},p,z^{1}%
,\ldots,z^{\underline{r}-1}$ is a new local coordinate system on $P$. In view
of Theorem \ref{TeorForger2}, $\omega$ is then multipresymplectic.

On the other hand, let $\omega$ be multipresymplectic. Then $\underline
{\omega}$ is polypresymplectic. Moreover, (around every point in $P$) there are
$\alpha$-adapted local coordinates%
\begin{equation}
\ldots,x^{i},\ldots,q^{A},\ldots,p_{A}^{i},\ldots,p,z^{1},\ldots,z^{r}
\label{Eq10}%
\end{equation}
such that, locally, $\omega=dp_{A}^{i}dq^{A}d^{n-1}x_{i}-dpd^{n}x$ and
$\underline{\omega}=d^{V}\!p_{A}^{i}d^{V}\!q^{A}\otimes d^{n-1}x_{i}$. This shows
that for all $y\in P$, 
\[
D_{y}^{\omega}=\langle\ldots,\left.  \tfrac{\partial
}{\partial z^{\alpha}}\right\vert _{y},\ldots\rangle\neq\underline{D}%
_{y}^{\omega}=\langle\ldots,\left.  \tfrac{\partial}{\partial z^{\alpha}%
}\right\vert _{y},\ldots,\left.  \tfrac{\partial}{\partial p}\right\vert
_{y}\rangle.
\]
 \end{proof}

\begin{remark}
\label{Rem1}Let $\omega$ be a PD prehamiltonian system. First of
all, notice that, if $\omega$ is a multipresymplectic structure then, in view
of Proposition \ref{Prop2}, PD Hamilton equations of $\omega$ do not
possess solutions. In this sense, multipresymplectic structures do not contain
any dynamical information.

Now, the proof of Proposition \ref{Prop2} also shows that if
$\underline{\omega}$ is a polypresymplectic structure and $\ker\omega
=\ker\underline{\omega}$ then $\omega$ is locally in the form
\[
\omega=dp_{A}^{i}dq^{A}d^{n-1}x_{i}-dHd^{n}x
\]
where $\tfrac{\partial H}{\partial z^{\alpha}}=0$, $\alpha=1,2,\ldots$, i.e.,
$H$ is constant along the leaves of the distribution $D^{\omega}=\underline
{D}^{\omega}$.
\end{remark}

\begin{example}
\label{Rem2}Let $\omega\in\Omega^{2}$ be a multisymplectic structure on
$\alpha$. In this case $\ker\omega=0$, while $\underline{D}^{\omega}$ is a
$1$-dimensional (involutive) distribution. Leaves of $\underline{D}^{\omega}$
are submanifolds in the fibers of $\alpha$. denote by $\underline{P}$ the set
of leaves of $\underline{D}^{\omega}$. There is an obvious projection
$\underline{\alpha}:\underline{P}\longrightarrow M$. Suppose that
$\underline{\alpha}:\underline{P}\longrightarrow M$ is a smooth fiber bundle
and $\mathfrak{p}:P\longrightarrow\underline{P}$ a smooth submersion (which is
always true locally). There is a distinguished class of (local) PD Hamiltonian
systems on $\underline{\alpha}$. Indeed, let $U\subset\underline{P}$ be an
open subbundle and let $\mathscr{H}:U\longrightarrow P$ be a local section of
$\mathfrak{p}$. Then $\omega^{\prime}:=\mathscr{H}^{\ast}(\omega)\in\Omega
^{2}(U,\underline{\alpha})$ is a PD Hamiltonian system. In particular, if we
choose coordinates on $P$ as in (\ref{Eq10}) (here $r=0$), then $\ldots
,x^{i},\ldots,q^{A},\ldots,p_{A}^{i},\ldots$ are coordinates on $\underline
{P}$, $\mathscr{H}$ is given by
\[
\mathscr{H}^{\ast}(p)=H,
\]
$H$ being a local function on $\underline{P}$, and $\omega^{\prime}$ is
locally given by
\[
\omega^{\prime}=dp_{A}^{i}dq^{A}d^{n-1}x_{i}-dHd^{n}x,
\]
in particular $\theta^{\prime}:=p_{A}^{i}dq^{A}d^{n-1}x_{i}-Hd^{n}x$ is a
local potential of $\omega^{\prime}$. Finally, PD Hamilton equations of
$\omega^{\prime}$ read
\begin{align*}
q_{i}^{A}  &  =\tfrac{\partial H}{\partial p_{A}^{i}},\\
p_{A}^{i}{}_{,i}  &  =-\tfrac{\partial H}{\partial q^{A}},
\end{align*}
which are de Donder-Weyl equations (see, for instance, \cite{gs73}).
\end{example}

\subsection{PD Hamiltonian Systems and Variational Calculus}

We show in this subsection that PD Hamilton equations are locally variational.
First of all, an element $\theta\in\Omega^{1}$ may be understood as a
(fiber-wise affine) horizontal $n$-form over $J^{1}\alpha$, i.e., as an
element $\mathscr{L}^{\theta}\in\overline{\Lambda}{}^{n}(J^{1}\alpha
,\alpha_{1})$ via
\[
\mathscr{L}_{c}^{\theta}:=i_{c}\theta_{y},\quad c\in J^{1}\alpha
,\;y=\alpha_{1,0}(c)\in P.
\]
In its turn $\mathscr{L}^{\theta}$ is a 1st order Lagrangian density in the
fiber bundle $\alpha$ determining an action functional which we denote by
$S^{\theta}=\int\mathscr{L}^{\theta}$. If $\theta$ is locally
given by $\theta=\theta_{a}^{i}dy^{a}d^{n-1}x_{i}-Hd^{n}x$, $\ldots,\theta_a^i,\ldots,H$ being local functions on $P$, then
$\mathscr{L}^{\theta}$ is locally given by $\mathscr{L}^{\theta}=L^{\theta
}d^{n}x$, where $L^{\theta}$ is the local function on $J^{1}\alpha$ given by
\[
L^{\theta}=(\theta_{b}^{i}y{}_{i}^{b}-H).
\]
In particular, if $\theta=\delta\nu$ for some $\nu\in\Omega^{0}=\overline
{\Lambda}{}^{n-1}$ locally given by $\nu=\nu^{i}d^{n-1}x_{i}$, then
\begin{equation}
L^{\theta}=(\partial_{i}+y_{i}^{a}\partial_{a})\nu^{i}, \label{Eq34}%
\end{equation}
i.e., $L^{\theta}$ is a total divergence.

\begin{prop}
Let $\omega\in\Omega^{2}$ be a $\delta$-exact PD prehamiltonian system. Then
 PD Hamilton equations of $\omega$ coincide with Euler-Lagrange equations
associated with the action $S^{\theta}:=\int\mathscr{L}^{\theta}$, where
$\theta\in\Omega^{1}$ is the opposite of any potential of $\omega$, i.e.,
$-\delta\theta=\omega$. Moreover, if $H^{1}(\Omega,\delta)=0$ then $S^{\theta
}$ is independent of the choice of $\theta$ and does only depend on $\omega$.
\end{prop}

\begin{proof}
The first part of the proposition can be proved in local coordinates. Indeed,
compute variational derivatives of $L^{\theta}$,
\begin{align*}
\tfrac{\delta}{\delta y^{b}}L^{\theta}  &  :={\partial}_{b}%
L^{\theta}-(\partial_{i}+y_{i}^{a}\partial_{a})\tfrac{\partial}{\partial
y{}_{i}^{b}}L^{\theta}\\
&  =y{}_{i}^{a}(\partial_{b}\theta_{a}^{i}-\partial_{a}\theta_{b}%
^{i})-\partial_{a}H-\partial_{i}\theta_{a}^{i}\\
&  =-2\omega_{ab}^{i}y{}_{i}^{a}+\omega_{b}.
\end{align*}
where we used (\ref{Eq32}). To prove the second part of the proposition, use
(\ref{Eq34}) to conclude that, for $\nu\in\Omega^{0}$, $\delta L^{\delta
\nu}/\delta y^{a}=0$.
 \end{proof}

\begin{remark}
Condition $H^{1}(\Omega,\delta)=0$ depends on the topology of the fiber bundle
$\alpha$. It is satisfied, for instance, if $H^{n}(M)=0$ and $H^{1}%
(\mathcal{F})=0$, $\mathcal{F}$ being, as above, the abstract fiber of
$\alpha$. Indeed, if $H^{1}(\mathcal{F})=0$ then $H^{1}(\underline{\Omega
},\underline{\delta})=0$ so that, the first part of the exact sequence
(\ref{Eq28}) reads
\[
\xymatrix@C=32pt{ 0 \ar[r] & H^{0}(\Omega,\delta) \ar[r] & \Lambda^{n-1}(M) \ar[r]^-{d_M^{n-1}} & \Lambda^{n}(M) \ar[r] &
H^{1}(\Omega,\delta) \ar[r] & 0}
\]
and $H^{1}(\Omega,\delta)\simeq H^{n}(M)=0$.
\end{remark}

\section{PD Noether Symmetries and Currents\label{SecNoether}}

\subsection{PD Noether Theorem and PD Poisson Bracket}

The multisymplectic analogues of Hamiltonian vector fields and Poisson
bracket in symplectic geometry have been longly investigated
\cite{ks75,k97,fr01,fpr03,fpr03b,fr05}. We here propose the natural
definitions for general PD Hamiltonian systems. Notice that, even if they look
formally identical to (or possibly less general than) the ones proposed in
\cite{ks75,fr01,fpr03,fpr03b}, our definitions have actually got a dynamical
content, not only a kinematical one (see Remark \ref{Rem1}), so that, for
instance, we can prove a PD version of (Hamiltonian) Noether theorem. That's why, e.g., we
will better speak about PD Noether symmetries rather than Hamiltonian
(multi)vector fields \cite{fpr05}.

Let $\omega$ be a PD prehamiltonian system on the bundle $\alpha
:P\longrightarrow M$. In the following we assume $\alpha$ to have connected fiber.

\begin{definition}\label{DefPDNoether}
Let $Y\in V\mathrm{D}$ and $f\in\Omega^{0}$. If $i_{Y}\omega=\delta f$, then
$Y$ and $f$ are said to be a \emph{PD Noether symmetry} and a \emph{PD Noether
current} of $\omega$ (relative to each other), respectively.
\end{definition}

Denote by $\mathscr{S}(\omega)$ and $\mathscr{C}(\omega)$ the sets of
PD Noether symmetries and PD Noether currents of $\omega$, respectively. A
PD Noether symmetry $Y$ (relative to a PD Noether current $f$) is a symmetry
of $\omega$ in the sense that
\[
L_{Y}\omega=i_{Y}\delta\omega+\delta i_{Y}\omega=\delta\delta f=0.
\]
The next proposition clarifies in what sense a PD Noether current is a
conserved current for $\omega$.

\begin{prop}
[PD -Noether theorem]Let $Y\in\mathscr{S}(\omega)$ and $f\in\mathscr{C}(\omega
)$ be a PD Noether symmetry and a PD Noether current of $\omega$ relative to
each other. Then $\sigma^{\ast}(f)\in\Lambda^{n-1}(M)$ is a closed form for
every solution $\sigma$ of PD Hamilton equations.
\end{prop}

\begin{proof}
First of all, let $\varrho\in\Omega^{1}$ and let $\tau$ be a (local) section of
$\alpha$. It is easy to show (for instance, using local coordinates) that
$\tau^{\ast}(\varrho)=i_{\dot{\tau}}\varrho|_{\tau}\in\Lambda^{n}(M)$. Then
\begin{align*}
d\sigma^{\ast}(f)  &  =\sigma^{\ast}(df)\\
&  =\sigma^{\ast}(\delta f)\\
&  =i_{\dot{\sigma}}\delta f|_{\sigma}\\
&  =i_{\dot{\sigma}}i_{Y}\omega|_{\sigma}\\
&  =i_{Y|_{\sigma}}i_{\dot{\sigma}}\omega|_{\sigma}\\
&  =0.
\end{align*}

 \end{proof}

We are now in the position to introduce a Lie bracket among PD Noether currents.

\begin{prop}
Let $Y_{1},Y_{2}\in\mathscr{S}(\omega)$ be PD Noether symmetries relative to
the PD Noether currents $f_{1},f_{2}\in\mathscr{C}(\omega)$, respectively.
Then $[Y_{1},\allowbreak Y_{2}]\in\mathscr{S}(\omega)$ and $f:=L_{Y_{1}}%
f_{2}\in\mathscr{C}(\omega)$ and they are relative to each other. Moreover,
$f$ is independent of the choice of $Y_{1}$ among the PD Noether symmetries
relative to the PD Noether current $f_{1}$.
\end{prop}

\begin{proof}
Compute
\begin{align*}
\delta L_{Y_{1}}f_{2}  &  =L_{Y_{1}}\delta f_{2}\\
&  =L_{Y_{1}}i_{Y_{2}}\omega\\
&  =i_{[Y_{1},Y_{2}]}\omega+i_{Y_{2}}L_{Y_{1}}\omega\\
&  =i_{[Y_{1},Y_{2}]}\omega.
\end{align*}
Now, let $V\in\ker\omega$. Then $L_{V}f_{2}=i_{V}\delta f_{2}=i_{V}i_{Y_{2}%
}\omega=0$. This proves the second part of the proposition.
 \end{proof}

Let $Y_{1},Y_{2},f_{1},f_{2}$ be as in the above proposition.

\begin{prop}
The $\mathbb{R}$-bilinear map 
\[
\mathscr{C}(\omega)\times\mathscr{C}(\omega
)\ni(f_{1},f_{2})\longmapsto\{f_{1},f_{2}\}:=L_{Y_{1}}f_{2}\in H(\omega),
\]
$Y_{1}$ being a PD Noether symmetry relative to $f_{1}$, is a Lie bracket.
\end{prop}

\begin{proof}
Let $Y_2 \in \mathscr{S}(\omega)$ be a PD Noether symmetry relative to $f_2 \in \mathscr{C}(\omega)$. Skew-symmetry of $\{\cdot,\cdot\}$ immediately follows from the remark:
\begin{align*}
\{f_{1},f_{2}\}  &  =L_{Y_{1}}f_{2}\\
&  =i_{Y_{1}}\delta f_{2}+\delta i_{Y_{1}}f_{2}\\
&  =i_{Y_{1}}i_{Y_{2}}\omega.
\end{align*}
Now, check Leibniz rule. Let $Y_{3}\in \mathscr{S}(\omega)$ and $f_3 \in \mathscr{C}(\omega)$ be another pair of
PD Noether symmetry, PD Noether current relative to each other. Then
\begin{align*}
\{f_{1},\{f_{2},f_{3}\}\}  &  =L_{Y_{1}}\{f_{2},f_{3}\}\\
&  =L_{Y_{1}}L_{Y_{2}}f_{3}\\
&  =L_{[Y_{1},Y_{2}]}f_{3}+L_{Y_{2}}L_{Y_{1}}f_{3}\\
&  =\{\{f_{1},f_{2}\},f_{3}\}+\{f_{2},\{f_{1},f_{3}\}\}.
\end{align*}
 \end{proof}

PD Noether symmetries and PD Noether currents of a PD Hamiltonian system
constitute very small Lie subalgebras of the Lie algebras of higher symmetries
and conservation laws of PD Hamilton equations, for which there have been given
fully satisfactory definitions and have been developed many infinite jet based
computational techniques \cite{b...99}. Nevertheless, it is worthy to give Definition \ref{DefPDNoether} and
to carefully analyze it, independently on infinite jets, in view of the possibility
of developing a \textquotedblleft (multi)symplectic theory\textquotedblright\ of
higher symmetries and conservation laws (see, for instance, \cite{v09}). In
Section \ref{SecComp} we propose some specific examples.

Finally, notice that, in general, nor a PD Noether current is uniquely
determined by the relative PD Noether symmetry nor vice versa (unless
$\ker\omega=0$). However, \textquotedblleft non-trivial PD Noether
symmetries\textquotedblright\ are in one to one correspondence with
\textquotedblleft non-trivial PD Noether currents\textquotedblright\ in the following sense. Clearly, $\ker
\omega\subset\mathscr{S}(\omega)$ and $H^{0}(\Omega,\delta)\subset
\mathscr{C}(\omega)$. We will call elements in $\ker\omega$ \emph{gauge
PD Noether symmetries} (see below) and elements in $H^{0}(\Omega,\delta)$
(i.e., closed $(n-1)$ -forms on $M$, see Corollary \ref{CorH0}) \emph{trivial
PD Noether currents}.

\begin{remark}
It is easy to see that $\ker\omega$ and $H^{0}(\Omega,\delta)$ are ideals in
the Lie algebras $\mathscr{S}(\omega)$ and $\mathscr{C}(\omega)$,
respectively. Let $\overline{\mathscr{S}}(\omega):=\mathscr{S}(\omega
)/\ker\omega$ and let $\overline{\mathscr{C}}(\omega):=\mathscr{C}(\omega
)/H^{0}(\Omega,\delta)$ be the quotient Lie algebras. Then the map
\[
\overline{\mathscr{S}}(\omega)\ni Y+\ker\omega\longmapsto f+H^{0}%
(\Omega,\delta)\in\overline{\mathscr{C}}(\omega),
\]
where $Y\in\mathscr{S}(\omega)$ and $f\in\mathscr{C}(\omega)$ are relative to
each other, is a well defined isomorphism of Lie algebras. It is natural to
call elements in $\overline{\mathscr{S}}(\omega)$ and $\overline
{\mathscr{C}}(\omega)$ \emph{non-trivial PD Noether symmetries} and
\emph{non-trivial PD Noether currents}, respectively. Indeed, elements in
$\ker\omega$ are trivial symmetries in that they are infinitesimal gauge
transformations (see next subsection), and elements in $H^{0}(\Omega,\delta)$ are
trivial conserved currents in that they are conserved currents for every
PD prehamiltonian system $\omega$, independently of $\omega$.
\end{remark}

\subsection{Gauge Reduction of PD Hamiltonian Systems}

From a physical point of view, elements in $\ker\omega$ are infinitesimal
gauge transformations and therefore should be quotiented out via a reduction
of the system. In this section we assume $\ker\omega=\ker\underline{\omega}$
or, which is the same, $\operatorname{Ker}\omega\neq\varnothing$. As a further
regularity condition we assume that the leaves of $D^{\omega}=\underline
{D}^{\omega}$ form a smooth fiber bundle $\widetilde{P}$ over $M$, whose
projection we denote by $\widetilde{\alpha}:\widetilde{P}\longrightarrow M$,
in such a way that the canonical projection $\mathfrak{p}:P\longrightarrow
\widetilde{P}$ is a smooth bundle. The last condition is always fulfilled at
least locally. Notice, also, that, by construction, $\mathfrak{p}$ has
connected fiber.

\begin{thrm}
There exists a unique PD Hamiltonian system $\widetilde{\omega}$ in
$\widetilde{\alpha}$ such that 1) $\omega=\mathfrak{p}{}^{\ast}(\widetilde
{\omega})$, 2) $\ker\widetilde{\omega}=\ker\underline{\widetilde{\omega}}=0$
and 3) a local section $\sigma$ of $\alpha$ is a solution of the PD Hamilton
equation of $\omega$ iff $\mathfrak{p}\circ\sigma$ (which is a local section
of $\widetilde{\alpha}$) is a solution of PD Hamilton equations of
$\widetilde{\omega}$.
\end{thrm}

\begin{proof}
Let $\widetilde{\nabla}\in C(\widetilde{P},\widetilde{\alpha})$. There exists
a (non-unique) connection $\nabla\in C(P,\alpha)$ such that $\nabla$ and
$\widetilde{\nabla}$ are $\mathfrak{p}$-compatible. To prove this, choose a
connection $\square$ in $\mathfrak{p}$ and lift the planes of $\widetilde
{\nabla}$ to $P$ by means of $\square$. It is easy to show that the so
obtained distribution on $P$ defines a connection $\nabla$ in $\alpha$ with
the required property. Similarly, every vector field $\widetilde{X}\in
V\mathrm{D}(\widetilde{P},\widetilde{\alpha})$ can be lifted to a (non-unique)
$\mathfrak{p}$-projectable vector field $X\in V\mathrm{D}(P,\alpha)$
such that $\widetilde{X}$ is its projection. Then $X\in\mathrm{D}%
_{V}(P,\mathfrak{p})$. Consider $\eta:=\omega(\nabla)(X)\in\Omega^{0}%
(P,\alpha)$ and prove that $L_{Y}\eta=0$ for any $Y\in V\mathrm{D}%
(P,\mathfrak{p})$. Indeed, let $Y\in V\mathrm{D}(P,\mathfrak{p})$. Then
$[Y,X]\in V\mathrm{D}(P,\mathfrak{p})$. Similarly $[\![Y,H_{\nabla}]\!]\in\overline{\Lambda}%
{}^{1}(P,\alpha)\otimes V\mathrm{D}(P,\mathfrak{p})\subset\overline{\Lambda}%
{}^{1}(P,\alpha)\otimes V\mathrm{D}(P,\alpha)$. Now, $V\mathrm{D}%
(P,\mathfrak{p})=\ker\omega$ by construction, and therefore
\begin{align*}
L_{Y}\eta &  =L_{Y}i_{X}i_{\nabla}\omega\\
&  =[L_{Y},i_{X}]i_{\nabla}\omega+i_{X}L_{Y}i_{\nabla}\omega\\
&  =i_{[Y,X]}i_{\nabla}\omega+i_{X}i_{\nabla}L_{Y}\omega+i_{X}[L_{Y}%
,i_{\nabla}]\omega\\
&  =i_{\nabla}i_{[Y,X]}\omega+i_{X}i_{[\![Y,H_{\nabla}]\!]}\omega\\
&  =0.
\end{align*}
Since fibers of $\mathfrak{p}$ are connected we conclude that $\eta
=\mathfrak{p}^{\ast}(\widetilde{\eta})$ for a unique $\widetilde{\eta}%
\in\Omega^{0}(\widetilde{P},\widetilde{\alpha})$. Put
\[
\widetilde{\omega}(\widetilde{\nabla})(\widetilde{X}):=\widetilde{\eta},
\]
so that $\widetilde{\omega}$ is a well defined element in $\Omega^{2}(\widetilde
{P},\widetilde{\alpha})$. Indeed, let $\nabla^{\prime}\in C(P,\alpha)$ be also
$\mathfrak{p}$-compatible with $\widetilde{\nabla}$ and let $X^{\prime}\in
V\mathrm{D}(P,\alpha)$ be another $\mathfrak{p}$-projectable vector field projecting onto $\widetilde{X}$. Then $\nabla^{\prime}-\nabla
\in\overline{\Lambda}{}^{1}(P,\alpha)\otimes V\mathrm{D}(P,\mathfrak{p})$ and
$X^{\prime}-X\in V\mathrm{D}(P,\mathfrak{p})$. Therefore,
\begin{align*}
\omega(\nabla^{\prime})(X^{\prime})  &  =i_{X^{\prime}}i_{\nabla^{\prime}%
}\omega\\
&  =i_{X^{\prime}}i_{\nabla}\omega+i_{X^{\prime}}i_{\nabla^{\prime}-\nabla
}\underline{\omega}\\
&  =i_{X}i_{\nabla}\omega+i_{X^{\prime}-X}i_{\nabla}\omega\\
&  =i_{X}i_{\nabla}\omega+i_{\nabla}i_{X^{\prime}-X}\omega\\
&  =\omega(\nabla)(X).
\end{align*}
Moreover, $\omega=\mathfrak{p}^{\ast}(\widetilde{\omega})$ by construction.

Let us compute $\ker\underline{\widetilde{\omega}}$. Thus, let $\widetilde
{X}\in V\mathrm{D}(\widetilde{P},\widetilde{\alpha})$ be such that
$i_{\widetilde{X}}\underline{\widetilde{\omega}}=0$ and let $X\in V\mathrm{D}%
(P,\alpha)$ be as above. Then $i_{X}\underline{\omega}=\mathfrak{p}^{\ast
}(i_{\widetilde{X}}\underline{\widetilde{\omega}})=0$. This shows that $X\in$
$V\mathrm{D}(P,\mathfrak{p})$ and then $\widetilde{X}=0$.

Finally, let $\sigma$ be a local section of $\alpha$, $\widetilde{\sigma
}:=\mathfrak{p}\circ\sigma$, $\widetilde{X}\in V\mathrm{D}(\widetilde
{P},\widetilde{\alpha})$ and let $X$ be as above. Compute
\begin{align*}
(i_{\widetilde{\sigma}{}^{\cdot}}\widetilde{\omega}|_{\widetilde{\sigma}%
})(\widetilde{X}|_{\widetilde{\sigma}})  &  =i_{\widetilde{\sigma}{}^{\cdot}%
}(i_{\widetilde{X}}\widetilde{\omega})|_{\widetilde{\sigma}}\\
&  =\widetilde{\sigma}{}^{\ast}(i_{\widetilde{X}}\widetilde{\omega})\\
&  =(\sigma^{\ast}\circ\mathfrak{p}^{\ast})(i_{\widetilde{X}}\widetilde
{\omega})\\
&  =\sigma^{\ast}(i_{X}\omega)\\
&  =i_{\dot{\sigma}}(i_{X}\omega)|_{\sigma}\\
&  =(i_{\dot{\sigma}}\omega|_{\sigma})(X|_{\sigma}).
\end{align*}
This shows that $i_{\dot{\sigma}}\omega|_{\sigma}=0$ iff $i_{\widetilde
{\sigma}{}^{\cdot}}\widetilde{\omega}|_{\widetilde{\sigma}}=0$.
 \end{proof}

\begin{prop}
There are natural isomorphisms of Lie algebras
\begin{align*}
\overline{\mathscr{S}}(\omega)  &  \simeq\mathscr{S}(\widetilde{\omega}),\\
\mathscr{C}(\omega)  &  \simeq\mathscr{C}(\widetilde{\omega}).
\end{align*}

\end{prop}

\begin{proof}
First of all let $f\in\mathscr{C}(\omega)$ and $X\in\mathscr{S}(\omega)$ be
relative to each other. Then $f=\mathfrak{p}^{\ast}(\widetilde{f})$ for some
$\widetilde{f}\in\Omega^{0}(\widetilde{P},\widetilde{\alpha})$ and $X$ is
$\mathfrak{p}$-projectable. Indeed, for all $Y\in\ker\omega$,
\[
L_{Y}f=i_{Y}\delta f+\delta i_{Y}f=i_{Y}i_{X}\omega=i_{[Y,X]}\omega=0.
\]
Moreover,
\[
\mathfrak{p}^{\ast}(\delta\widetilde{f})=\delta\mathfrak{p}^{\ast}%
(\widetilde{f})=\delta f=i_{X}\omega=\mathfrak{p}^{\ast}(i_{\widetilde{X}%
}\widetilde{\omega}),
\]
where $\widetilde{X}$ denotes the $\mathfrak{p}$-projection of $X$, and,
therefore, $\delta\widetilde{f}=i_{\widetilde{X}}\widetilde{\omega}$, i.e.,
$\widetilde{f}\in\mathscr{C}(\widetilde{\omega})$ and $\widetilde{X}%
\in\mathscr{S}(\widetilde{\omega})$ is a PD Noether symmetry relative to it.
Thus, maps
\begin{align}
\overline{\mathscr{S}}(\omega)\ni X+\ker\omega &  \longmapsto\widetilde{X}%
\in\mathscr{S}(\widetilde{\omega}),\label{Eq36}\\
\mathscr{C}(\omega)\ni f  &  \longmapsto\widetilde{f}\in\mathscr{C}(\widetilde
{\omega}). \label{Eq37}%
\end{align}
are well defined. Conversely, let $\widetilde{X}_{1}\in\mathscr{S}(\widetilde{\omega})$,
$\widetilde{f}_{1}\in\mathscr{C}(\widetilde{\omega})$ be relative to each
other, $X_{1}\in V\mathrm{D}(P,\alpha)$ be any $\mathfrak{p}$-projectable
vector field, $\widetilde{X}_{1}\in V\mathrm{D}(\widetilde{P},\widetilde
{\alpha})$ be its projection, and let $f_{1}:=\mathfrak{p}^{\ast
}(\widetilde{f}_{1})\in\Omega^{0}(P,\alpha)$. Then $X_{1}\in\mathscr{S}(\omega
)$ and $f_{1}\in\mathscr{C}(\omega)$ is a PD Noether current relative to it.
Indeed,
\[
i_{X_{1}}\omega=\mathfrak{p}^{\ast}(i_{\widetilde{X}_{1}}\widetilde{\omega
})=\mathfrak{p}^{\ast}(\delta\widetilde{f}_{1})=\delta\mathfrak{p}^{\ast
}(\widetilde{f}_{1})=\delta f_{1}.
\]
We conclude that (\ref{Eq36}) and (\ref{Eq37}) are inverted by
\begin{align*}
\mathscr{S}(\widetilde{\omega})\ni\widetilde{X}_{1}  &  \longmapsto X_{1}%
+\ker\omega\in\overline{\mathscr{S}}(\omega),\\
\mathscr{C}(\widetilde{\omega})\ni\widetilde{f}_{1}  &  \longmapsto f_{1}%
\in\mathscr{C}(\omega),
\end{align*}
respectively.
 \end{proof}

\section{Examples\label{SecComp}}

\subsection{Non-Degenerate Examples}

Let $\alpha:\mathbb{R}^{2n+1}\ni(x^{1},\ldots,x^{n},u,u_{1},\ldots
,u_{n})\longmapsto(x^{1},\ldots,x^{n})\in\mathbb{R}^{n}$, $n>1$. Consider
$T,V\in C^{\infty}(\mathbb{R}^{2n+1})$ of the form $T=T(u_{1},\ldots,u_{n})$
and $V=V(u)$, respectively. The form
\[
\omega:=\tfrac{\partial^{2}T}{\partial u_{i}\partial u_{j}}du_{i}%
(dud^{n-1}x_{j}-u_{j}d^{n}x)-V^{\prime}dud^{n}x,
\]
is a PD prehamiltonian system on $\alpha$ (here and in what follows a prime
\textquotedblleft$\;{}^{\prime}\;$\textquotedblright\ denotes differentiation
with respect to $u$). The associated PD Hamilton equations read
\[%
\begin{array}
[c]{c}%
\tfrac{\partial^{2}T}{\partial u_{i}\partial u_{j}}\partial_{j}u_{i}%
+V^{\prime}=0,\\
\partial_{i}u=u_{i},
\end{array}
\]
which are in turn equivalent to%
\begin{align}
\tfrac{\partial^{2}T}{\partial u_{i}\partial u_{j}}\partial_{ij}%
^{2}u+V^{\prime}  &  =0,\label{Eq14}\\
\partial_{i}u  &  =u_{i},\nonumber
\end{align}
$\partial_{ij}^{2} := \partial_i \partial_j$, $i,j=1,\ldots,n$.
Moreover, for%
\[
\det\left(  \tfrac{\partial^{2}T}{\partial u_{i}\partial u_{j}}\right)
\neq0,
\]
$\omega$ is a PD Hamiltonian system. We will only consider this case in the
following. Thus, put $T^{ij}:=\tfrac{\partial^{2}T}{\partial u_{i}\partial
u_{j}}$, $i,j=1,\ldots,n$, and let $(T_{ij})$ be the inverse matrix of
$\left(  T^{ij}\right)  $. As examples notice that

\begin{enumerate}
\item For $T=\tfrac{1}{2}g^{ij}u_{i}u_{j}$,
\[
(g^{ij})=\left(
\begin{array}
[c]{cccc}%
-1 & 0 & \cdots & 0 \\
0 & 1 & \cdots & 0\\
\vdots & \vdots & \ddots & \vdots\\
0 & 0 & \cdots & 1
\end{array}
\right)  ,
\]
(resp., $g^{ij}=\delta^{ij}$, $i,j=1,\ldots,n$), (\ref{Eq14}) reduces to the
wave equation (resp., the Poisson equation) with a $u$-dependent potential $V$
(including the $f$-Gordon equation as a particular example, if $n=2$ and
$f=-V^{\prime}$).

\item For $n=2$, $T=\sqrt{1+g^{ij}u_{i}u_{j}}$, $g^{ij}=\delta^{ij}$,
$i,j=1,2$, and $V=0$, (\ref{Eq14}) reduces to the equation for minimal
surfaces in $\mathbb{R}^{3}$ transversal to the projection $\mathbb{R}^{3}%
\ni(x_{1},x_{2},u)\longmapsto(x_{1},x_{2})\in\mathbb{R}^{2}$.
\end{enumerate}

Let us search for PD Noether symmetries and currents of $\omega$. Let
$Y=U\tfrac{\partial}{\partial u}+U_{i}\tfrac{\partial}{\partial u^{i}}%
\in\mathrm{V}D$ and let $f=f^{i}d^{n-1}x_{i}\in\Omega^{0}$. Then
\begin{align*}
i_{Y}\omega &  =T^{ij}(U_{i}du-Udu_{i})d^{n-1}x_{j}-\left(  T^{ij}u_{i}%
U_{j}+V^{\prime}U\right)  d^{n}x,\\
\delta f  &  =\partial_{i}f^{i}d^{n}x+\tfrac{\partial}{\partial u}%
f^{i}dud^{n-1}x_{i}+\tfrac{\partial}{\partial u_{k}}f^{i}du_{k}d^{n-1}x_{i}.
\end{align*}
Recall that $Y$ and $f$ are a PD Noether symmetry and a PD Noether current
relative to each other, respectively, iff $i_{Y}\omega=\delta f$, i.e.,
\begin{align}
\partial_{i}f^{i}+T^{ij}u_{i}U_{j}+V^{\prime}U  &  =0,\label{Eq16}\\
\tfrac{\partial}{\partial u}f^{i}-T^{ij}U_{j}  &  =0,\label{Eq17}\\
\tfrac{\partial}{\partial u_{j}}f^{i}+T^{ij}U  &  =0. \label{Eq18}%
\end{align}
It follows from (\ref{Eq18}) that $\tfrac{\partial}{\partial u_{j}}%
f^{i}=\tfrac{\partial}{\partial u_{i}}f^{j}$, $i,j=1,\ldots,n$, and then
\[
\tfrac{\partial^{2}}{\partial u_{k}\partial u_{j}}f^{i}=\tfrac{\partial^{2}%
}{\partial u_{k}\partial u_{i}}f^{j},\quad i,j,k=1,\ldots,n.
\]
Now,
\begin{align*}
\tfrac{\partial^{2}}{\partial u_{k}\partial u_{j}}f^{i}  &  =\tfrac{\partial
}{\partial u_{k}}\tfrac{\partial}{\partial u_{j}}f^{i}\\
&  =-\tfrac{\partial}{\partial u_{k}}\left(  T^{ij}U\right) \\
&  =-\tfrac{\partial^{3}T}{\partial u_{k}\partial u_{i}\partial u_{j}}%
U-T^{ij}\tfrac{\partial}{\partial u_{k}}U
\end{align*}
Similarly,
\[
\tfrac{\partial^{2}}{\partial u_{k}\partial u_{i}}f^{j}=\tfrac{\partial
}{\partial u_{i}}\tfrac{\partial}{\partial u_{k}}f^{j}=-\tfrac{\partial^{3}%
T}{\partial u_{i}\partial u_{j}\partial u_{k}}U-T^{jk}\tfrac{\partial
}{\partial u_{i}}U.
\]
Therefore,
\[
T^{ij}\tfrac{\partial}{\partial u_{k}}U-T^{jk}\tfrac{\partial}{\partial u_{i}%
}U=0
\]
Contracting with $T_{ij}$ we find $(n-1)\tfrac{\partial}{\partial u_{i}}U=0$
and, therefore,
\[
U=U(x^{1},\ldots,x^{n},u),
\]
so that (\ref{Eq18}) can be rewritten as
\[
\tfrac{\partial}{\partial u_{j}}\left(  f^{i}+\tfrac{\partial T}{\partial
u_{i}}U\right)  =0.
\]
We conclude that
\begin{equation}
f^{i}=-\tfrac{\partial T}{\partial u_{i}}U+A^{i} \label{Eq19}%
\end{equation}
for some $A^{i}=A^{i}(x^{1},\ldots,x^{n},u)$, $i=1,\ldots,n$. Notice that
(\ref{Eq17}) can be used to determine the $U_{j}$'s from the $f^{i}$'s via
\[
U_{j}=T_{ji}\tfrac{\partial}{\partial u}f^{i}.
\]
It remains to solve (\ref{Eq16}) which, in view of (\ref{Eq19}), reduces to
\begin{equation}
\left(  \partial_{i}+u_{i}\tfrac{\partial}{\partial u}\right)  A^{i}%
-\tfrac{\partial T}{\partial u_{i}}\left(  \partial_{i}+u_{i}\tfrac{\partial
}{\partial u}\right)  U+V^{\prime}U=0. \label{Eq20}%
\end{equation}
We cannot go further on in solving (\ref{Eq20}) without better specifying $T$.
In the following we will only consider two special cases.

\begin{enumerate}
\item $T=\tfrac{1}{2}g^{ij}u_{i}u_{j}$, $(g^{ij})$ being a constant, non
degenerate, symmetric matrix with inverse $(g_{ij})$. In this case
(\ref{Eq20}) reads
\begin{equation}
\partial_{i}A^{i}+V^{\prime}U+\left(  \tfrac{\partial}{\partial u}A^{i}%
-g^{ij}\partial_{j}U\right)  u_{i}-\left(  g^{ij}\tfrac{\partial}{\partial
u}U\right)  u_{i}u_{j}=0. \label{Eq21}%
\end{equation}
The left hand side of (\ref{Eq21}) is polynomial in $u_{1},\ldots,u_{n}$.
Thus, all the corresponding coefficients must vanish, i.e.,
\begin{align}
\tfrac{\partial}{\partial u}U  &  =0,\label{Eq22}\\
\tfrac{\partial}{\partial u}A^{i}-g^{ij}\partial_{j}U  &  =0,\label{Eq23}\\
\partial_{i}A^{i}+V^{\prime}U  &  =0. \label{Eq24}%
\end{align}
From (\ref{Eq22}), $U=U(x^{1},\ldots,x^{n})$ and, then, from (\ref{Eq23}),
$\tfrac{\partial^{2}}{\partial u^{2}}A^{i}=0$, $i=1,\ldots,n$, which in turn
implies, using (\ref{Eq23}) again,
\[
A^{i}=\left(  g^{ij}\partial_{j}U\right)  u+B^{i}%
\]
for some $B^{i}=B^{i}(x^{1},\ldots,x^{n})$. Finally, (\ref{Eq24}) implies
\[
\left(  g^{ij}\partial_{ij}^{2}U\right)  u+\partial_{i}B^{i}+V^{\prime}U=0
\]
and differentiating once more with respect to $u$
\[
g^{ij}\partial_{ij}^{2}U+V^{\prime\prime}U=0.
\]
Since $U$ doesn't depend on $u$, if

\begin{enumerate}
\item $V^{\prime\prime\prime}\neq0$. Then $U=0$ so that
\[
f^{i}=\tfrac{1}{2}\partial_{j}B^{ji},\quad U_{j}=0
\]
for some $B^{ij}=-B^{ji}=B^{ij}(x^{1},\ldots,x^{n})$, i.e.,
\[
Y=0\quad\text{and}\quad f=d\beta,
\]
$\beta=B^{ji}d^{n-2}x_{ji}$, where $d^{n-2}x_{ji}:=i_{\partial_{j}}%
d^{n-1}x_{i}$, $i,j=1,\ldots,n$. Therefore, $\omega$ doesn't posses
PD Noether symmetries nor non-trivial PD Noether currents.

\item $V^{\prime\prime\prime}=0$. Then $V=\tfrac{1}{2}\mu u^{2}$ for some
constant $\mu$ and
\[
g^{ij}\partial_{ij}^{2}U+\mu U=0,\quad f^{i}=g^{ij}(u\,\partial_{j}%
U-u_{j}U)+\tfrac{1}{2}\partial_{j}B^{ji},\quad U_{j}=\partial_{j}U,
\]
for some $B^{ij}=-B^{ji}=B^{ij}(x^{1},\ldots,x^{n})$. Thus,
\[
Y=U\tfrac{\partial}{\partial u}+\partial_{j}U\tfrac{\partial}{\partial u_{j}%
}\quad\text{and}\quad f=g^{ij}(u\,\partial_{j}U-u_{j}U)d^{n-1}x_{i}+d\beta,
\]
$\beta=B^{ji}d^{n-2}x_{ji}$, where $U$ is any solution of the PD Hamilton
equation
\begin{equation}
g^{ij}\partial_{ij}^{2}u+\mu u=0. \label{Eq25}%
\end{equation}
Let us compute the PD Poisson bracket. Consider two solutions of (\ref{Eq25}),
say $U_{1},U_{2}$, the corresponding PD Noether symmetries $Y_{1},Y_{2}$ and
associated PD Noether currents $f_{1},f_{2}$. Then
\[
\{f_{1},f_{2}\}=L_{Y_{1}}f_{2}=g^{ij}(U_{1}\partial_{j}U_{2}-U_{2}\partial
_{j}U_{1})d^{n-1}x_{i},
\]
which, as can be easily checked, is a trivial conservation law.
\end{enumerate}

\item $n=2$, $T=\sqrt{1+\delta^{ij}u_{i}u_{j}}$ and $V=0$. In this case
(\ref{Eq20}) reads
\begin{equation}
\tau^{1/2}\left(  \partial_{i}+u_{i}\tfrac{\partial}{\partial u}\right)
A^{i}=\delta^{ij}u_{j}\left(  \partial_{i}+u_{i}\tfrac{\partial}{\partial
u}\right)  U, \label{Eq26}%
\end{equation}
where $\tau=1+\delta^{ij}u_{i}u_{j}$. Squaring both sides of (\ref{Eq26}) we
get
\[
\tau\left[  \left(  \partial_{i}+u_{i}\tfrac{\partial}{\partial u}\right)
A^{i}\right]  ^{2}-\left[  \delta^{ij}u_{j}\left(  \partial_{i}+u_{i}%
\tfrac{\partial}{\partial u}\right)  U\right]  ^{2}=0,
\]
whose left hand side is polynomial in $u_{1},u_{2}$. Collecting homogeneous
terms we get
\begin{align}
&  \left[  \left(  \tfrac{\partial}{\partial u}U\right)  ^{2}\delta
^{ij}-\left(  \tfrac{\partial}{\partial u}A^{i}\right)  \left(  \tfrac
{\partial}{\partial u}A^{j}\right)  \right]  \delta^{kl}u_{i}u_{j}u_{k}%
u_{l}\nonumber\\
+  &  2\delta^{ij}\left[  \delta^{kl}\left(  \tfrac{\partial}{\partial
u}U\right)  \left(  \partial_{l}U\right)  -\left(  \tfrac{\partial}{\partial
u}A^{k}\right)  \left(  \partial_{l}A^{l}\right)  \right]  u_{i}u_{j}%
u_{k}\nonumber\\
-  &  \left[  \delta^{ij}\left(  \partial_{k}A^{k}\right)  ^{2}+\left(
\tfrac{\partial}{\partial u}A^{i}\right)  \left(  \tfrac{\partial}{\partial
u}A^{j}\right)  -\delta^{ik}\delta^{jl}\left(  \partial_{k}U\right)  \left(
\partial_{l}U\right)  \right]  u_{i}u_{j}\nonumber\\
+  &  2\left(  \partial_{j}A^{j}\right)  ^{2}\left(  \tfrac{\partial}{\partial
u}A^{i}\right)  u_{i}+\left(  \partial_{i}A^{i}\right)  ^{2}=0. \label{Eq27}%
\end{align}
All coefficient of the left hand side of (\ref{Eq27}) must vanish. It follows
that
\[
\tfrac{\partial}{\partial u}U=\partial_{1}U=\partial_{2}U=0,\quad
\tfrac{\partial}{\partial u}A^{1}=\tfrac{\partial}{\partial u}A^{2}%
=0,\quad\partial_{1}A^{1}+\partial_{2}A^{2}=0,
\]
i.e., $U$ is a constant while $A^{1}=\partial_{2}B$, $A^{2}=-\partial_{1}B$
for some $B=B(x^{1},x^{2})$. Thus,
\[
Y=U\tfrac{\partial}{\partial u},\quad f=U\,\tau^{-1/2}\left(  u_{2}%
dx^{1}-u_{1}dx^{2}\right)  +dB.
\]
It is obvious that the PD Poisson bracket is also trivial in this case.
\end{enumerate}

\subsection{A Degenerate, Constrained Example}

The example in this subsection is taken from \cite{gmr09}. Let $\alpha:\mathbb{R}^{3m+2}\times\mathbb{R}_{+}\ni(q^{1},\ldots,q^{m}%
,s_{1},\ldots,s_{m},t_{1},\ldots,t_{m},s,t;e)\longmapsto(s,t)\in\mathbb{R}%
^{2}$. The form
\[
\omega:=-dt_{\alpha}dq^{\alpha}ds+ds_{\alpha}dq^{\alpha}dt-\delta^{\alpha
\beta}(et_{\alpha}dt_{\beta}-s_{\alpha}ds_{\beta})dsdt-\varepsilon dedsdt,
\]
where $\varepsilon:=\tfrac{1}{2}(\delta^{\alpha\beta}t_{\alpha}t_{\beta}-1)$,
is a PD prehamiltonian system on $\alpha$. The associated PD Hamilton
equation reads
\[%
\begin{array}
[c]{c}%
\tfrac{\partial}{\partial t}t_{\alpha}+\tfrac{\partial}{\partial s}s_{\alpha
}=0,\\
\tfrac{\partial}{\partial t}q^{\alpha}=e\delta^{\alpha\beta}t_{\beta},\\
\tfrac{\partial}{\partial s}q^{\alpha}=-\delta^{\alpha\beta}s_{\beta},\\
\varepsilon=0,
\end{array}
\]
$\alpha=1,\ldots,m$, which is in turn equivalent to
\[%
\begin{array}
[c]{c}%
e^{-1}\tfrac{\partial^{2}}{\partial t^{2}}q^{\alpha}-\tfrac{\partial^{2}%
}{\partial s^{2}}q^{\alpha}=e^{-2}\left(  \tfrac{\partial}{\partial
t}q^{\alpha}\right)  \left(  \tfrac{\partial}{\partial t}e\right)  ,\\
e^{2}=\delta_{\alpha\beta}\left(  \tfrac{\partial}{\partial t}q^{\alpha
}\right)  \left(  \tfrac{\partial}{\partial t}q^{\beta}\right)  ,\\
t_{\alpha}=e^{-1}\delta_{\alpha\beta}\tfrac{\partial}{\partial t}q^{\beta},\\
s_{\alpha}=\delta_{\alpha\beta}\tfrac{\partial}{\partial s}q^{\beta},
\end{array}
,
\]
Notice that $\underline{D}^{\omega}$ is generated by $\tfrac{\partial
}{\partial e}$, while
\[
D_{y}^{\omega}=\left\{
\begin{array}
[c]{cc}%
\mathbf{0} & \text{for }\varepsilon(y)\neq0\\
\left\langle \left.  \tfrac{\partial}{\partial e}\right\vert _{y}\right\rangle
& \text{for }\varepsilon(y)=0
\end{array}
\right.  ,\quad y\in P
\]
we conclude that $P_{(1)}$ is the hypersurface defined by $\delta^{\alpha\beta
}t_{\alpha}t_{\beta}=1$. It is easy to see that, actually, $\breve{P}=P_{(1)}$.

Let us search for PD Noether symmetries and currents of $\omega$. Let
$Y=Q^{\alpha}\tfrac{\partial}{\partial q^{\alpha}}+S_{\alpha}\tfrac{\partial
}{\partial s_{\alpha}}+T_{\alpha}\tfrac{\partial}{\partial t_{\alpha}}%
+E\tfrac{\partial}{\partial e}\in V\mathrm{D}$ and let $f=\alpha ds+\beta
dt\in\Omega^{0}$. Then $i_{Y}\omega=\delta f$ iff
\begin{equation}%
\begin{array}
[c]{c}%
\tfrac{\partial}{\partial s}\beta-\tfrac{\partial}{\partial t}\alpha
=\delta^{\beta\gamma}(s_{\beta}S_{\gamma}-et_{\beta}T_{\gamma})-\varepsilon
E,\\
\tfrac{\partial}{\partial q^{\alpha}}\alpha=-T_{\alpha},\quad\tfrac{\partial
}{\partial q^{\alpha}}\beta=S_{\alpha},\quad\tfrac{\partial}{\partial
t_{\alpha}}\alpha=\tfrac{\partial}{\partial s_{\alpha}}\beta=Q^{\alpha}\\
\tfrac{\partial}{\partial s_{\alpha}}\alpha=\tfrac{\partial}{\partial
t_{\alpha}}\beta=0,\quad\tfrac{\partial}{\partial e}\alpha=\tfrac{\partial
}{\partial e}\beta=0,
\end{array}
\label{Eq35}%
\end{equation}
$\alpha=1,\ldots,m$. Equations (\ref{Eq35}) can be easily solved and give
quite large $\mathscr{S}(\omega)$ and $\mathscr{C}(\omega)$. Namely,
\[
\alpha=C^{\alpha}t_{\alpha}+A,\quad\beta=C^{\alpha}s_{\alpha}+B,
\]
and
\[%
\begin{array}
[c]{c}%
Q^{\alpha}=C^{\alpha},\quad T_{\alpha}=-\tfrac{\partial C^{\beta}}{\partial
q^{\alpha}}t_{\beta}-\tfrac{\partial A}{\partial q^{\alpha}},\quad S_{\alpha
}=\tfrac{\partial C^{\beta}}{\partial q^{\alpha}}s_{\beta}+\tfrac{\partial
B}{\partial q^{\alpha}},\\
\varepsilon E=\tfrac{\partial C^{\alpha}}{\partial s}s_{\alpha}-\tfrac
{\partial C^{\alpha}}{\partial t}t_{\alpha}+\tfrac{\partial B}{\partial
s}-\tfrac{\partial A}{\partial t}-\delta^{\alpha\beta}\left[  s_{\alpha
}\left(  \tfrac{\partial C^{\gamma}}{\partial q^{\beta}}s_{\gamma}%
+\tfrac{\partial B}{\partial q^{\beta}}\right)  +et_{\alpha}\left(
\tfrac{\partial C^{\gamma}}{\partial q^{\beta}}t_{\gamma}+\tfrac{\partial
A}{\partial q^{\beta}}\right)  \right]
\end{array}
\]
where $A,B,\ldots,C^{\alpha},\ldots,D^{\alpha\beta},\ldots,E^{\alpha},\ldots$
are arbitrary functions of the only $s,t,\dots,q^{\beta},\ldots$.

Compute the PD Poisson bracket. Let $f_{1}$,$f_{2}$ be PD Noether currents
determined by functions $A_{1},B_{1},\ldots,C_{1}^{\alpha},\ldots$ and
$A_{2},B_{2},\ldots,C_{2}^{\alpha},\ldots$ respectively. A straightforward
computation shows that
\[
\{f_{1},f_{2}\}=(C^{\alpha}t_{\alpha}+A)ds+(C^{\alpha}s_{\alpha}+B)dt
\]
with
\begin{align*}
A  &  =C_{1}^{\beta}\tfrac{\partial}{\partial q^{\beta}}A_{2}-C_{2}^{\beta
}\tfrac{\partial}{\partial q^{\beta}}A_{1},\\
B  &  =C_{1}^{\beta}\tfrac{\partial}{\partial q^{\beta}}B_{2}-C_{2}^{\beta
}\tfrac{\partial}{\partial q^{\beta}}B_{1},\\
C^{\alpha}  &  =C_{1}^{\beta}\tfrac{\partial}{\partial q^{\beta}}C_{2}%
^{\alpha}-C_{2}^{\beta}\tfrac{\partial}{\partial q^{\beta}}C_{1}^{\alpha},
\end{align*}
$\alpha=1,\ldots,m$.

\subsection{A Degenerate, Unconstrained Example}

Finally, we propose an example of reduction. Consider the cotangent bundle
$\pi:T^{\ast}\mathbb{M}\ni A_{i}dx^{i}|_{(x^{1},\ldots,x^{n})}\longmapsto
(x^{1},\ldots,x^{n})\in\mathbb{M}$ and let $\alpha:=\pi_{1}:(x^{1},\dots
,x^{n},\ldots,A_{i},\ldots,A_{i,j},\ldots)\ni J^{1}\pi\longmapsto(x^{1}%
,\dots,x^{n})\in\mathbb{M}$, $\mathbb{M}$ being the $n$-dimensional Minkowski
space. As such $\mathbb{M}$ is endowed with the metric $g:=g_{ij}dx^{i}\cdot
dx^{j}$ where
\[
(g_{ij})=\left(
\begin{array}
[c]{cccc}%
-1 & 0 & \cdots & 0\\
0 & 1 & \cdots & 0\\
\vdots & \vdots & \ddots & \vdots\\
0 & 0 & \cdots & 1
\end{array}
\right)  .
\]
In the following we will raise and lower indexes using $g$. Let
\[
\omega:=2dA^{[j,i]}\left(  \tfrac{1}{2}A_{i,j}d^{n}x-dA_{i}d^{n-1}%
x_{j}\right)
\]
Then $\omega$ is a PD prehamiltonian system on $\pi$ whose PD Hamilton
equation reads
\begin{align*}
\partial_{k}A^{[i,k]}  &  =0,\\
\partial_{\lbrack j}A_{i]}  &  =A_{[i,j]},
\end{align*}
$i,j=1,\ldots,n$, which are equivalent to Maxwell equations for the vector
potential
\begin{align*}
(\partial_{k}\partial^{k})A_{i}-\partial_{i}\partial_{k}A^{k}  &  =0,\\
A_{[i,j]}  &  =\partial_{\lbrack j}A_{i]}.
\end{align*}
Notice that
\[
\ker\omega=\ker\underline{\omega}=\left\langle \ldots,\tfrac{\partial
}{\partial A_{i,j}}+\tfrac{\partial}{\partial A_{j,i}},\ldots\right\rangle .
\]
Therefore $\omega$ is \textquotedblleft degenerate and
unconstrained\textquotedblright. Moreover, leaves of $D^{\omega}=\underline
{D}^{\omega}$ are given by $A_{[i,j]}=\mathrm{const}$. We conclude that
$J^{1}\pi$ \textquotedblleft reduces\textquotedblright\ via
\begin{align*}
\mathfrak{p}:J^{1}\pi &  \longrightarrow T^{\ast}\mathbb{M\times}_{M}%
\wedge^{2}T^{\ast}\mathbb{M\simeq R}^{n(n+3)/2}\\
(x^{1},\dots,x^{n},\ldots,A_{i},\ldots,A_{i,j},\ldots)  &  \longmapsto
(x^{1},\dots,x^{n},\ldots,A_{i},\ldots,F_{ij},\ldots)
\end{align*}
where $F_{ij}=F_{[ij]}$, $\mathfrak{p}^{\ast}(F_{ij}):=2A_{[j,i]}$ and
$\omega=\mathfrak{p}^{\ast}(\widetilde{\omega})$, with
\[
\widetilde{\omega}:=dF^{ij}\left(  \tfrac{1}{4}F_{ji}d^{n}x-dA_{i}d^{n-1}%
x_{j}\right)
\]
is a PD Hamiltonian system on
\[
\widetilde{\alpha}:\mathbb{R}^{n(n+3)/2}\ni(x^{1},\dots,x^{n},\ldots
,A_{i},\ldots,F_{ij},\ldots)\longmapsto(x^{1},\dots,x^{n})\in\mathbb{R}^{n},
\]
whose PD Hamilton equations read
\begin{align*}
\partial_{k}F^{ik}  &  =0,\\
\partial_{\lbrack j}A_{i]}  &  =2F_{ji},
\end{align*}
which are Maxwell equations for the field strength.

\section{PD Hamiltonian Systems in Mathematical Physics and Geometry}\label{SecMPG}

A system of PDEs is multisymplectic if it is in the form
\begin{equation}
K_{ab}^{i}\partial_{i}y^{a}=\partial_{b}H, \label{msPDE}%
\end{equation}
with $K_{ab}^{i}=-K_{ba}^{i}(y^{1},\ldots,y^{m})$ and $H=H(y^{1},\ldots
,y^{m})$ given functions and $\kappa^{i}:=K_{ab}^{i}dy^{a} dy^{b}$ a
symplectic form for all $i=1,\ldots,n$. Multisymplectic systems of PDEs where
first introduced in \cite{b97} to study the interaction and stability of
non-linear waves. More generally, a multi-symplectic formulation of a PDE
proved to be a useful tool in the stability analysis. They have been also
introduced multisymplectic integrators \cite{br01}, generalizing the
symplectic methods so popular in numerical Hamiltonian dynamics.

Notice that multisymplectic PDEs (\ref{msPDE}) are actually PD Hamilton
equations of the PD Hamiltonian system
\begin{equation}
\omega=\dfrac{1}{2}\kappa^{i} d^{n-1}x_{i}+dH d^{n}x.
\label{msPDHS}%
\end{equation}
PD\ Hamiltonian systems of the form (\ref{msPDHS}) are of a special kind.
Indeed, they are autonomous in two respects. Neither their kinematics (encoded
by the symplectic forms $\kappa^{i}$) nor their dynamics (encoded by the
\textquotedblleft Hamiltonian function\textquotedblright\ $H$) depend on
space-time coordinates. While the latter condition can be always achieved by
adding auxiliary coordinates, and, therefore, it is not a lack of generality,
the former is a special feature of multisymplectic PDEs among general PD
Hamilton equations. We conclude that PD Hamiltonian systems provide an
intrinsic (coordinate independent) geometric formalization of the theory of
Multisymplectic PDEs.

Many equations of fluid dynamics (and, more generally, of continuum
mechanics), including the Euler equation as an instance, are multisymplectic
PDEs, and, therefore, PD\ Hamilton equations (possibly after a suitable change
of coordinates) \cite{m...01,chh07}. Notice that fluid dynamics on a general
Riemannian manifold may not possess a multisymplectic formulation, while it
still possesses a PD Hamiltonian one, since, in this case, the
\textquotedblleft kinematics\textquotedblright\ depends on the metric and,
therefore, on the space-time.

Systems of hydrodynamic type, with their Dubrovin-Novikov Poisson structures
\cite{dn83,dn84}, are also multisymplectic (see, for instance, \cite{m98}),
and, therefore, PD Hamiltonian. For the latter systems, the relation between
the multisymplectic structure and the Poisson structure has been discussed at
least in the integrable, 1-dimensional case of the KdV equation (see
\cite{g88} for details).

As already mentioned, a Lagrangian field theory in the bundle $\pi
:E\longrightarrow M$, with Lagrangian density $\mathscr{L}$ locally given by
$\mathscr{L}=Ld^{n}x$, $L=L(x^{1},\ldots,x^{n},\ldots,u_{i}^{\alpha},\ldots)$,
where the $u_{i}^{\alpha}$'s denotes \textquotedblleft space--time
derivatives\textquotedblright\ of the field variables $u^{\alpha}$'s,
determines canonically a PD Hamiltonian system $\omega_{\mathscr{L}}$ in the
bundle $J^{1}\pi\longrightarrow M$, locally given by
\[
\omega_{\mathscr{L}}=d\tfrac{\partial L}{\partial u_{i}^{\alpha}}
du^{\alpha} d^{n-1}x_{i}-dE d^{n}x,\quad E:=u_{i}^{\alpha}%
\tfrac{\partial L}{\partial u_{i}^{\alpha}}-L.
\]
In many cases the PD Hamilton equations of $\omega_{\mathscr{L}}$ are
equivalent to the Euler-Lagrange equations (even in presence of gauge
symmetries). The (functional) space of solutions of the Euler-Lagrange
equations carries a canonical (pre)symplectic structure $\boldsymbol{\omega}$
\cite{z87,cw87} whose degeneracy distribution is made of gauge symmetries
\cite{lw90}. Therefore, the (classical) gauge reduction of a field theory
basically amounts to the symplectic reduction of $\boldsymbol{\omega}$. This
remark is at the very basis of the BV formalism for the quantization of gauge
theories \cite{ht92}. Interesting examples may be found at the frontier
between Theoretical Physics and Geometry. For instance, the reduction of the
space of flat connections in a principal bundle over a Riemannian 3-manifold
to the \emph{moduli space} of (gauge equivalent) flat connections amounts to
the symplectic reduction of the presymplectic structure of the Chern-Simons
theory. Similarly, the Atiyah-Hitchin manifold is a moduli space of (gauge
equivalent) magnetic monopoles. As such it inherits the symplectic structure
from the presymplectic structure of the Yang-Mills-Higgs theory. This is
precisely the symplectic sector of the Atiyah-Hitchin hyper-K\"{a}hler
structure. Now, $\boldsymbol{\omega}$ can be understood as a cohomology class
of a suitable complex and $\omega_{\mathscr{L}}$ as a cocycle representing it
\cite{v09}. For completeness, we present the (reduced) PD Hamiltonian system
corresponding to the Yang-Mills-Higgs theory.

Let $\mathbb{M}$ be the $n$-dimensional Minkowski space with coordinates
$\ldots,x^{i},\ldots$ and let $\mathfrak{g}$ be a Lie algebra with a basis
$\ldots,\sigma_{a},\ldots$. Consider the following (vector) bundles:

\begin{enumerate}
\item $\mathbb{M}\times$ $\mathfrak{g}\longrightarrow$ $\mathbb{M}$, with
linear bundle coordinates $\ldots,\phi^{a},\ldots$;

\item a copy of $T^{\ast}\mathbb{M}\otimes\mathfrak{g}\longrightarrow$
$\mathbb{M}$, with bundle coordinates $\ldots,\psi_{i}^{a},\ldots$;

\item a second copy of $T^{\ast}\mathbb{M}\otimes\mathfrak{g}\longrightarrow$
$\mathbb{M}$, with bundle coordinates $\ldots,A_{i}^{a},\ldots$;

\item $\wedge^{2}T^{\ast}\mathbb{M}\otimes\mathfrak{g}\longrightarrow$
$\mathbb{M}$, with bundle coordinates $\ldots,F_{ij}^{a},\ldots$, $i<j$.

\item The fibered product $\alpha:P\longrightarrow M$ of the above.
\end{enumerate}

In $\alpha$ there is the following natural PD Hamiltonian system
\[
\omega_{YMH}=k_{ab}[dF^{a}{}^{ij}(dA_{j}^{b}+\tfrac{1}{2}[A_{l}%
,A_{j}]^{b}dx^{l})+d\psi^{a}{}^{i}(\phi^{b}+[A_{l},\phi]^{b}%
dx^{l})] d^{n-1}x_{i}-Ld^{n}x,
\]
where
\[
L=k_{ab}(\tfrac{1}{4}F^{a}{}^{ij}F_{ij}^{b}+\tfrac{1}{2}\psi^{a}{}^{i}\psi
_{i}^{b}),
\]
and $(k_{ab})$ is the Killing form of $\mathfrak{g}$. The PD Hamilton
equations of $\omega_{YMH}$ are
\begin{align*}
\partial_{i}\psi^{a}{}^{i}+[A_{i},\psi]^{a}  &  =0,\\
\partial_{i}\phi^{a}+[A_{i},\phi]^{a}  &  =\psi_{i}^{a},\\
\partial_{i}F^{a}{}^{ij}+[A_{i},F^{ij}]^{a}  &  =0,\\
\partial_{i}A_{j}^{a}-\partial_{j}A_{i}^{a}+[A_{i},A_{j}]^{a}  &  =F_{ij}^{a},
\end{align*}
which are basically the Yang-Mills-Higgs equations.

\end{document}